\documentclass{amsart}
\usepackage[english]{babel}
\usepackage{amsmath}
\usepackage{graphicx}
\usepackage[foot]{amsaddr}
\usepackage{xcolor}
\usepackage{subcaption}
\usepackage{float}
\usepackage{multirow}
\newcommand{\fig}[1]{Fig.~\ref{#1}}{\color{blue}}

\usepackage{comment}
\usepackage{stackengine,graphicx}
\usepackage{tikz}
\usepackage{overpic}
\usepackage{subcaption}
\usepackage[top=1in, bottom=1.25in, left=1.25in, right=1.25in]{geometry}
\usepackage{hyperref}
    \hypersetup{
        colorlinks   = true,
        allcolors   =blue    
    }
\usepackage{tikz}
\usetikzlibrary{calc,fit,backgrounds}
\usetikzlibrary{shapes.geometric, arrows}
\usetikzlibrary{chains,
                fit,
                positioning,
                shapes}
\tikzstyle{startstop} = [rectangle, rounded corners, minimum width=10cm, minimum height=1cm,text centered, draw=black, fill=red!30]
\tikzstyle{io} = [trapezium, trapezium left angle=70, trapezium right angle=110, minimum width=2cm, minimum height=1cm, text centered, text width=2cm, draw=black, fill=blue!30]
\tikzstyle{io1} = [trapezium, trapezium left angle=70, trapezium right angle=110, minimum width=3cm, minimum height=1cm, text centered, text width=3cm, draw=black, fill=green!30]
\tikzstyle{process} = [rectangle, minimum width=1cm, minimum height=1cm, text centered, text width=2.5cm, draw=black, fill=orange!30]
\tikzstyle{decision} = [diamond, minimum width=3cm, minimum height=0.2cm, text centered,text width=1.2cm, draw=black, fill=green!30]
\tikzstyle{decision1} = [diamond, minimum width=2.8cm, minimum height=0.1cm,text width=1.8cm, draw=black, fill=green!30]
\tikzstyle{arrow} = [thick,->,>=stealth]
\tikzstyle{matheq} = [node distance=8.75cm, text width=21em, minimum width=1.5cm, minimum height=2em, text centered]
\tikzstyle{startstop1} = [rectangle, rounded corners, minimum width=5.5cm, minimum height=1cm,text centered, draw=black, fill=red!30]
\tikzstyle{io2} = [trapezium, trapezium left angle=70, trapezium right angle=110, minimum width=2cm, minimum height=1cm, text centered, text width=2cm, draw=black, fill=blue!30]
\tikzstyle{io1} = [trapezium, trapezium left angle=70, trapezium right angle=110, minimum width=3cm, minimum height=1cm, text centered, text width=3cm, draw=black, fill=green!30]
\tikzstyle{process1} = [rectangle, minimum width=1cm, minimum height=1.5cm, text centered, text width=2.5cm, draw=black, fill=orange!30]

\def\u{{\bm u}}

\def\U{{\bm U}}

\def\0{\boldsymbol{0}}

\newcommand{\bm}[1]{\mbox{\boldmath{$#1$}}}

\begin{document}

\title[]{A LSTM-enhanced surrogate model to simulate the dynamics of particle-laden fluid systems}
\author{Arash Hajisharifi$^1$, Rahul Halder$^1$, Michele Girfoglio$^1$, Andrea Beccari$^2$, Domenico Bonanni$^2$ and Gianluigi Rozza$^1$}
\address{$^1$ mathLab, Mathematics Area, SISSA, via Bonomea 265, I-34136 Trieste, Italy}
\address{$^2$ Dompé Farmaceutici SpA, EXSCALATE, Via Tommaso De Amicis, 95, I-80131 Napoli, Italy}

\begin{abstract}

The numerical treatment of fluid-particle systems is a very challenging problem because of the complex coupling phenomena occurring between the two phases. Although an accurate mathematical modelling is available to address this kind of applications, the computational cost of the numerical simulations is very expensive. The use of the most modern high performance computing infrastructures could help to mitigate such an issue but not completely to fix it. In this work we develop a non intrusive data-driven reduced order model (ROM) for Computational Fluid Dynamics (CFD) - Discrete Element Method (DEM)
simulations. The ROM is built using the proper orthogonal decomposition (POD) for the computation of the reduced basis space and the Long Short-Term Memory (LSTM) network for the computation of the reduced coefficients. We are interested to deal both with system identification and prediction. 
The most relevant novelties rely on (i) a filtering procedure of the full order snapshots to reduce the dimensionality of the reduced problem and (ii) a preliminary treatment of the particle phase. 
The accuracy of our ROM approach is assessed against the classic Goldschmidt fluidized bed benchmark problem.  
Finally we also provide some insights about the efficiency of our ROM approach.
\end{abstract}

\maketitle

\textbf{Keywords}: CFD-DEM, proper orthogonal decomposition, reduced order model, long short-term memory, data-driven techniques, fluidized bed.  

\section{Introduction}

The comprehension of the physical phenomena arising in fluid-solid systems is crucial in several chemical engineering processes. Experimental measurements are very hard to be obtained because of high temperature and pressure, as well as the presence of toxic substances. In addition, the cost of the instruments able to provide accurate measurements is rather expensive. 

For such reasons, numerical simulations represent a very useful tool to corroborate the information provided by the experimental tests and to support the design of industrial plants. Within this framework a particular importance is covered by the so-called Reduced Order Models (ROMs) \cite{hesthaven2016certified, quarteroni2015reduced, benner2015survey, benner2017model, benner2021system, bader, rozza2022advanced})
which are able to reduce the computational cost of the numerical simulations without a
significant loss of accuracy. The idea on which ROM is based is that the system at hand exhibits an intrinsic dimension which is much lower than the number of degrees of freedom resulting by the discretization of the system itself. The basic ROM framework consists of two steps. The first one is the so-called \emph{offline} phase, where a database of several solutions is collected by solving the original governing equations system (the so-called Full Order Model (FOM)) for
selected parameter values. The second step is the \emph{online} phase, during which the information
obtained in the offline phase is used to quickly compute the solution for newly specified values of the parameters. 
ROMs consist of two main different kinds: intrusive and non-intrusive. This is strictly related to the nature of the solver that is adopted. In the former case, one has access to the source
code of the FOM solver so that it is possible to manipulate directly the equations describing the
original problem. In the latter one, one can work only on the solutions, without any information about the original set of equations. In general, non-intrusive ROMs represent the ideal choice when one deals with closed source solvers that are widely used for industrial applications. Another advantage is that non-intrusive ROMs are able to provide higher computational speed-up especially when dealing with complex systems exhibiting a non-linear behaviour. See, e.g., \cite{hajisharifi2023non, hajisharifi2024comparison, siena2023data, rahman2019nonintrusive, salavatidezfouli2023modal, sheidani2023assessment, bakhshaei2024optimized}. 


In this work we deal with the development of a non intrusive data-driven ROM for Computational Fluid Dynamics (CFD) - Discrete Element Method (DEM) approach. It is an Eulerian-Lagrangian technique used for the simulation of
systems involving the interaction between a fluid flow and solid particles \cite{noro, hajisharifi2021particle, hajisharifi2022interface, review, clayton, heat}. The ROM which we propose is based on a Proper Orthogonal Decomposition (POD) - Long-Short Time Memory (LSTM) technique where the POD is used to extract the reduced basis space while the LSTM network is adopted for the computation of the reduced coefficients. LSTM has been used a lot for system prediction in several applications fields \cite{maulik2021reduced, pawar2020data, li2019deep, wang2018model, salavatidezfouli2024predictive}. We deal with the investigation of both system identification and prediction. 

At the best of our knowledge, there is not a great quantity of papers in literature about the development of ROMs for CFD-DEM simulations. A preliminary analysis based on a POD-Galerkin approach has been presented in \cite{palacios, yuan2005reduced}. 
Next in \cite{shuo, shuo2} 
the authors have proposed a non intrusive data-driven ROM framework based on the so-called PODI (Proper Orthogonal Decomposition with Interpolation) approach both for Eulerian and Lagrangian variables where the interpolation procedure for the evaluation of the reduced coefficients is based on Radial Basis Functions (RBFs). A further extension of such a work has been realized in \cite{hajisharifi2023non} where we have introduced a sensitivity analysis of the ROM error at varying of the number of POD modes retained. In addition, we also have performed a physical parametrization with respect to the Stokes number. Finally in \cite{li2023data} it has been performed a comparison between POD and Dynamic Mode Decomposition (DMD). All these works are related to the identification of the time evolution of the system. In \cite{li2023physics} the authors have proposed a Physics-informed DMD approach for system prediction. Anyway, the method exhibits several drawbacks: it has been validated exclusively for the Eulerian phase and provides reliable results only for a narrow time frame. 
When we were finalizing this manuscript, we became aware of \cite{li2024reduced} (appeared on March 15, 2024) where the authors compared PODI and POD-LSTM techniques both for identification and prediction of the Eulerian phase. The content of \cite{li2024reduced} shows strong analogies with this manuscript, but we could not be aware about that during the preparation of the present work. Since the adoption of LSTM for
fluid-particle systems is an innovative line of research, the results in \cite{li2024reduced} are not exhaustive. 
Here is how we complement them in this work: 
\begin{itemize}


\item A pre-processing step of the snapshots based on a Fast Fourier Transform (FFT) approach to reduce the dimensionality of the reduced problem deleting frequency content that would otherwise
be difficult or impossible to capture with the reduced model without a significant loss of accuracy and helping LSTM for better prediction. A similar procedure, but based on a Gaussian filtering approach, was introduced in \cite{farcas2022filtering} within the development of a data-driven ROM for rotating detonation engine simulations. 
\item A preliminary investigation about the prediction of the Lagrangian phase. At the best of our knowledge, this is the first paper dealing with such an aspect. All the works currently available concerning ROM prediction for CFD-DEM numerical simulations are limited to the study of the Eulerian phase. So we underline the difficulty to deal with the Lagrangian phase in the context of CFD-DEM
modelling that represents an open challenge.

\end{itemize}

The work is organized as follows. In Sec. \ref{sec:FOM} a brief description of the CFD-DEM model (i.e. our FOM) is reported. 
Then in Sec. \ref{sec:ROM} the ROM approach is described and the achieved results are introduced and discussed in Sec. \ref{sec:res}. Finally in Sec. \ref{sec:conc} conclusions are drawn and future perspectives are envisioned.


\section{The full order model}\label{sec:FOM}

Within the CFD-DEM technique, the fluid and solid phases are modeled based on Eulerian and Lagrangian approach, respectively, where the coupling between phases is given by the fluid-particle interaction forces \cite{Xu, tsu}. 

The governing equations of the fluid phase are outlined in Sec.~\ref{sec:fluid} while Sec.~\ref{sec:solid} shows the governing equations of the solid phase, including some insights about the corresponding numerical discretization. For further details the reader is referred to \cite{hajisharifi2023non}.


\subsection{Governing equations for the fluid-phase flow}\label{sec:fluid}

The volume-averaged Navier-Stokes equations (NSE) read as \cite{clayton}:

\begin{equation}
\centering
\frac{\partial \epsilon}{\partial t} + \nabla \cdot ( \epsilon \u ) = 0 \hspace{6cm} \text{in} \ \Omega \ \times \ (t_0, T],
\label{continuity}
\end{equation}
\begin{equation}
\frac{\partial (\epsilon \u)}{\partial t} + \nabla \cdot ( \epsilon \u \otimes \u) = - \nabla P - S_p + \nabla \cdot (\epsilon \bm{\tau}) + \epsilon \textbf{g} \hspace{1cm} \text{in} \ \Omega \ \times  \ (t_0, T],
\label{momentu}
\end{equation}
where $\Omega$ a is a fixed spatial domain and $(t_0, T]$ is the time interval of interest. In addition, $\epsilon$ is the fluid volume fraction, $\u$ is the fluid velocity, $\bm g$ is the gravity and $P = p/\rho_f$ is the modified pressure with $p$ and $\rho_f$ being the pressure and fluid density. The viscous stress tensor $\bm{\tau}$ is computed as follows: 

\begin{equation}
\bm{\tau} = \bm{\tau}_1 + \bm{\tau}_2 = \nu_f \left( \nabla\u + \nabla\u^T \right)- \dfrac{2}{3} \nu_f \left(\nabla \cdot \u \right) \bm{I},
\label{tau}
\end{equation}
where $\nu_f$ = $\mu_f / \rho_f$ represent the kinematic viscosity and $\bm{I}$ is the identity matrix.

The computational domain $\Omega$ is partitioned into $N_c$ cells or control volumes $\Omega_i$ with $i = 1, \dots, N_c$. 
The fluid volume fraction $\epsilon_{i}$, indicating the portion of the cell $i$ occupied by fluid, is defined as \cite{weller}: 

\begin{equation}
\epsilon_{i} = 1 - \frac {\sum_{j=1}^{n_p} \widetilde{\Omega}_j}{\Omega_{i}},
\label{eps}
\end{equation}
where  $n_p$ is the number of particles in the cell $i$ and $\widetilde{\Omega}_j$ represents the volume of the particle $j$. The interaction between fluid and particles is modelled by the source term $S_p$ which in the  cell $i$ is computed as \cite{zhu2007discrete}: 

\begin{equation}
S_{p,i} = \frac{\sum_{j=1}^{n_p} ( \mathbf{F}_{d,j} + \mathbf{F}_{\nabla p,j})}{\rho_f{\Omega}_{i}}.
\label{source}
\end{equation}

Here, $\mathbf{F}_{d,j}$ and $\mathbf{F}_{\nabla p,j}$ are the drag and pressure gradient forces, respectively \cite{zhu2007discrete}. 




The time interval of interest $(t_0, T]$ is partitioned using the time step $\Delta t \in \mathbb{R}^+$, resulting in time levels $t^n = t_0 + n \Delta t$, where $n$ ranges from 0 to the total number of time steps $N_{T}$. Eq.~\eqref{continuity} is discretized in time employing a first-order Euler scheme. On the other hand, all terms in eq.~\eqref{momentu} are treated implicitly except for the divergence term, which was discretized using a semi-implicit scheme. 

For the space discretization a second-order finite volume (FV) scheme is used. 
To address the pressure-velocity coupling, a partition approach is adopted, specifically the PIMPLE algorithm \cite{pimple}. This algorithm, which results by the combination of SIMPLE \cite{simple} and PISO \cite{issa} algorithms, offers a computationally robust framework. For the implementation of the numerical scheme we chose the finite volume C++ library OpenFOAM\textsuperscript{\textregistered} \cite{weller}. 


\subsection{Governing equations for the solid particle }\label{sec:solid}

The DEM model 
resolves particle motion 
by the translation and rotation second Newton's laws which are described as follows \cite{cundall1979discrete}:

\begin{equation}
m_j \frac{d\widetilde{\u}_j}{dt} =  \sum_{m=1}^{n_j^c} \bm{F}_{jm}^c +  \bm{F}_j^f + m_j \bm g ,\\
\label{translation}
\end{equation}
\begin{equation}
I_j \frac{d\mathbf \omega_j}{dt} = \sum_{m=1}^{n_j^c} \bm{M}_{jm}^c, 
\label{rotation}
\end{equation}
where 

\begin{equation}
m_j = \rho_p\dfrac{\pi d_p^3}{6} \quad \text{and} \quad I_j = \dfrac{m_j d_p^3}{6}
\end{equation}
are the mass and the moment of inertia of the  particle $j$, respectively. In addition, $\rho_p$,  $d_p$ and $\mathbf \omega_j$ are the density, the diameter and the angular velocity of the particle $j$, respectively.  $\bm{F}_{jm}^c$ and $\bm{M}_{jm}^c$ are the contact force and torque acting on particle $j$ by its $m$ contacts, which may involve other particles or wall \cite{fernandes}. The total number of contacts for the particle $j$ is denoted as $n_j^c$ and $\bm{F}_j^f = \bm{F}_{d,j} + \bm{F}_{\nabla p, j}$ is the particle-fluid interaction force acting on particle $j$. Notice that non-contact forces are not considered in this work.

The Stokes number, denoted as $Stk$, which characterizes the behavior of particles suspended in a fluid flow, is computed as \cite{ hajisharifi2022interface, marchioli2002mechanisms, sheidani2022study}:
\begin{equation}
    Stk = \frac{\tau_p}{\tau_f},
\label{Stk}
\end{equation}
where $\tau_f$ is the carrier fluid characteristic time and $\tau_p$ is the particle relaxation time  \cite{ hajisharifi2022interface, marchioli2002mechanisms, sheidani2022study}.
Eqs. \eqref{translation}-\eqref{rotation} are discretized by adopting a first-order Euler scheme. 


\section{The reduced order model}\label{sec:ROM}

We assume that any Eulerian (Lagrangian) variable can be approximated as a linear combination of a certain number of basis functions depending on the space $\bm{x}$ (label $l$ identifying the particles) only, multiplied by scalar coefficients that depend on the time and/or parameters of the problem at hand which can be physical or geometrical.

In this work the time $t$ is the only parameter of interest. In particular, we are going to for both
identification and prediction of the system dynamics. 
For what concerns the Eulerian phase, we are interested in the time evolution of the fluid volume fraction $\epsilon$ while for the Lagrangian phase we consider the position $\widetilde{\bm{x}}$ of the particles. Hence, the variables $(\epsilon, \widetilde{\bm{x}})$ are approximated by the reduced ones $(\epsilon_r, \widetilde{\bm{x}}_r)$ as follows 

\begin{align}
\epsilon \approx \epsilon_r = \sum_{i=1}^{N_{\epsilon}^r} \alpha_i(t) {\varphi}_i(\bm{x}), \quad 
\bm{\widetilde{x}} \approx \bm{\widetilde{x}}_r = \sum_{i=1}^{N_{\widetilde{x}}^r} \beta_i(t) \bm{\xi}_i(l).
\label{eq:podcoeff}
\end{align}

In eq. \eqref{eq:podcoeff}, $\bm{\alpha}(t) = [\alpha_i(t)]_{i=1}^{N_{\epsilon}^r}$ and $\bm{\beta}(t) = [\beta_i(t)]_{i=1}^{N_{\widetilde{x}}^r}$ are the coefficients of the reduced solutions, ${\varphi}$ and $\bm{\xi}$ are the basis functions and $N_{\epsilon}^r$ and $N_{\widetilde{x}}^r$ denote the cardinality of the reduced basis. In this work, we employ the Proper Orthogonal Decomposition with Long Short-Term Memory network (POD-LSTM) approach \cite{AHalder2020,ahmed2021nudged,tabib2022hybrid,deng2019time} which has been widely adopted for the development of ROMs in several engineering fields, e.g. aeroelastic and hydrodynamic applications \cite{AHalder2020,halder2023deep}, turbulent flow control \cite{mohan2018deep,rahman2019nonintrusive} and biomedical applications \cite{lyu2017long}.  


The POD-LSTM approach is based on the following \emph{offline}-\emph{online} paradigm: 

\begin{itemize}
\item{during the \emph{offline} phase, 
the time-dependent full order solutions belonging to a high dimensional manifold  $\mathcal{M}$ are collected for a given range of temporal instants and a global reduced basis for the space of the reduced solutions $\mathcal{M}_{POD}$ is extracted via POD. Then the full order snapshots are projected onto the POD space in order to compute the reduced coefficients and the LSTM training is performed to learn them. Due to the large number of degrees of freedom, the offline phase is computationally expensive. However, it is carried out only once.}

\item {during the \emph{online} phase the reduced coefficients for future time instances are quickly obtained from the learnt LSTM model. The reduced solution is then computed based on Eq. \eqref{eq:podcoeff}.} \end{itemize}


Next, we are going to describe the building blocks of the above algorithm.

\subsection{Proper Orthogonal Decomposition}
Let $\Phi (t^i)$ with $i = 1, \dots, N_t$ and $\Phi = \{\epsilon, \widetilde{\bm{x}}\}$ be the full order solutions obtained for different time instants $t^i$. Then we collect the snapshots in the matrix $\mathbb{\bm{S}} = [\Phi(t^1), \dots, \Phi(t^{N_t})]$. The application of the singular value decomposition to the matrix $\mathbb{\bm{S}}$ provides:\\
    \begin{equation}
    \mathbb{\bm{S}} = \U \bm{\Sigma} \bm{V}^T,
    \end{equation}
    where $\U \in \mathbb R^{{N_c}\times {N_t}}$ and $\bm{V} \in R^{{N_t}\times {N_t}}$ are the matrices composed by the left singular vectors  and right singular vectors, respectively, while $\bm{\Sigma} \in \mathbb R ^{{N_c}\times {N_t}}$ is the diagonal matrix containing the eigenvalues $\sigma_i$. The
POD space is constructed by retaining the first $N_\Phi^r << N_t $ columns of matrix $\U$. 
The value of $N_\Phi^r$ is commonly chosen to meet a user-provided threshold $\delta$ for the cumulative
energy of the eigenvalues 
    \begin{equation}
    E = \frac{\sum_{i=1}^{N_\Phi^r}\sigma_i}{\sum_{i=1}^{N_t}\sigma_i} \ge \delta.
    \label{k_equation}
    \end{equation}

After constructing the POD space we can approximate the input snapshots by using eq. \eqref{eq:podcoeff}:\\
    \begin{equation}
    \Phi(t^i) \approx \sum_{L=1}^{N_\Phi^r}\chi_L (t^i) \phi_L, \quad \text{with} \quad i = 1, \dots, N_t, \quad 
    \end{equation}
    where the modal coefficients $\chi_L (t^i)$ 
    are the elements of the matrix $\bm{C} = \U_{N_{\Phi}^r}^T \bm{S} \in \mathbb{R}^{N_\Phi^r \times N_t}$. 
 
\subsection{ Long Short-Term Memory network}\label{subsec:NN}
To demonstrate the network architecture associated with the LSTM model, we first introduce a simple Artificial Neural Network (ANN), followed by its extension to a Recurrent Neural Network (RNN) and finally to an LSTM network.  

ANN is a modelling technique able to formulate a nonlinear functional relationship between input and output data: see, e.g., 
\cite{pichi2023artificial,pier,san2019artificial,pawar2019deep,hesthaven2018non}.  
In our case the input-output data are the pairs $(t_i,  \bm{\chi}(t^i))$ 
where $\bm{\chi}(t^i) = [\chi_L(t^i)]_{L=1}^{N_\phi^r}$
with $i = 1, \dots, N_t$. 
The input-output mapping for a single-layer feed-forward ANN can be written as: 
\begin{equation}\label{eq:ANN1}
\begin{aligned}
{\bm{\chi}(t^i)} = f_{\text{act}}(\bm{W} t^i + \bm{b}).
\end{aligned}
\end{equation}
In eq. \eqref{eq:ANN1} $f_{act}$ is the nonlinear activation function, $\bm{W}$ is the weight matrix and
$\bm{b}$ is the bias vector. If there are multiple layers, a given hidden layer $h_l$ is fed to the next layer $h_{l+1}$: 
\begin{equation}\label{eq:ANN2}
\begin{aligned}
h_{l+1} = f_{\text{act}}(\bm{W} h_l + \bm{b}),
\end{aligned}
\end{equation}
where $h_l$ is generated from the output of the previous layers: see, e.g., \cite{pawar2019deep}.
The training algorithm involves an optimization problem to set the parameters of the network, $\bm{W}$ and $\bm{b}$, by minimizing a loss function which is a distance metric
between the truth and predicted data. 

\begin{figure}[ht] \label{Fig1}
\centering
{\label{fig:Fig1a}\includegraphics[width=.90\linewidth]{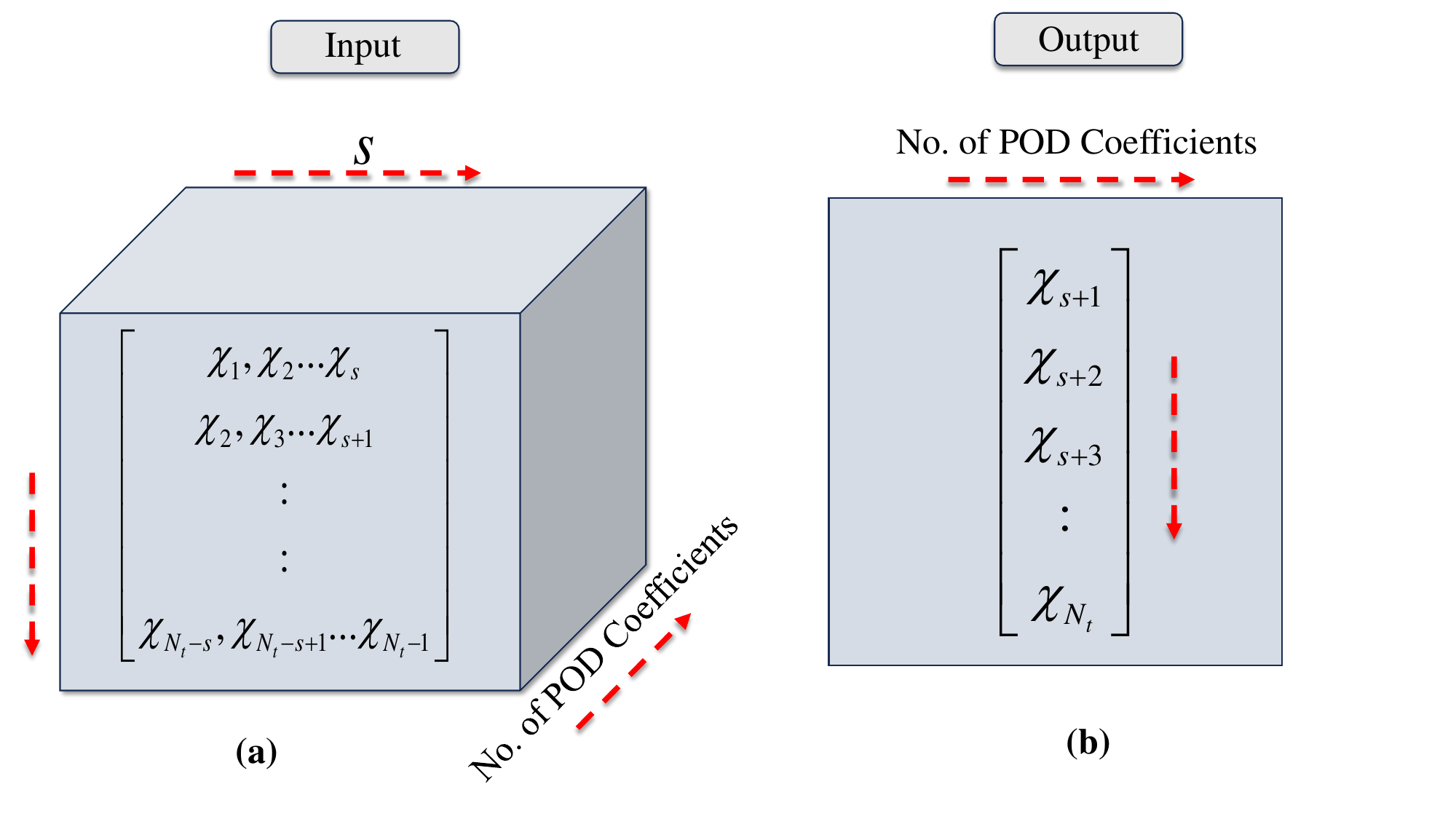}}
\caption{ 
Input (a) and output (b) RNN architecture.} 
\label{fig:inpoutRNN}
\end{figure}

Now, to model time prediction, we need to feed data from previous time steps back, i.e. from the time history of the system,  into the network. These networks with recurrent connections are called RNNs. 
Then we have 
\begin{equation}\label{eq:RNN1}
\begin{aligned}
h(t^i) = f_\text{act}\left(\bm{W} \bm{\chi}(t^{i}) +\bm{A} h(t^{i-1})+\bm{b}\right),
\end{aligned}
\end{equation}
 where $h(t^{i-1})$ is a hidden state that stores the sequence information within a certain number $s$ of previous data and $\bm{A}$ is an additional matrix which takes care of the feedback loop in the RNN architecture. 
The structure of the network is shown in Fig. \ref{fig:inpoutRNN}. For input data, it consists of a three-dimensional tensor: the first dimension 
is given by $N_t - s$, 
the second one by $s$ and the third one by $N_\Phi^r$. 
For output data, we have a two-dimensional tensor whose dimensions are $N_t-s+1$ and $N_\Phi^r$.


RNN layer creates instances of vanishing or exploding loss function gradient problems because of the long-term dependencies over time-sequential dataset. Then, the LSTM network is introduced in the RNN architecture to
reduce such an issue \cite{hochreiter1997long}. 
The LSTM introduces four distinct operations in each RNN cell: memory cell, input gate, forget gate, and output gate. For further details, the reader can refer to, e.g., \cite{AHalder2020,halder2023deep}.

In this work, the loss function associated with the LSTM network is given by the Mean Squared Error (MSE)
\begin{align}
    MSE = \frac{1}{N_t} \sum_{i=1}^{N_t} \left(\bm{\chi}(t^i) - \tilde{\bm{\chi}}(t^i) \right)^2,
\end{align}
where $\tilde{\bm{\chi}}(t^i)$ is the actual output vector
of the network. The ADAM optimization algorithm, which is based on a stochastic gradient
descent approach \cite{kingma2014adam}, is adopted for the
minimization of the loss function. 

After training the LSTM network we can compute system dynamics for a future time instance $t^{pred}$ by using eq. \eqref{eq:podcoeff}:\\
    \begin{equation}
    \Phi(t^{pred}) \approx \sum_{L=1}^{N_\Phi^r}\chi_{L} (t^{pred}) \phi_L. 
    \end{equation}

\subsection{Filtering of the snapshots}
Based on our previous work \cite{hajisharifi2023non}, a significant number of POD basis is required to describe accurately the dynamics of a CFD-DEM system. This could badly affect the efficiency of the ROM. 
To tackle this limitation and to reduce the computational complexity of our problem, before applying the POD to the high-fidelity snapshots, a Fast Fourier Transform (FFT) \cite{nussbaumer1982fast} algorithm is
employed to dump the high-frequency oscillations which increase the dimension of the latent space. 
Subsequently, the filtered data, rather than the original ones, will be used to build the ROM. As it is shown in Sec. \ref{sec:res}, by means of this methodology, we can maintain the same energy threshold with a fewer number of modes with increasing computational efficiency and without compromising the quality of the solution.

The application of the FFT to the fluid volume fraction $\epsilon$ and particle position $\widetilde{\bm{x}}$ provides

\begin{equation}\label{eq:DFT}
\begin{aligned}
\hat{\epsilon}_k= \text{FFT}(\epsilon) = \sum_{n=0}^{N_t-1} {\epsilon}_n^\star \cdot e^{-i 2 \pi \frac{k}{N_t} n}, \quad 
\hat{\bm{x}}_k= \text{FFT}(\bm{x}) = \sum_{n=0}^{N_t-1} {\widetilde{\bm{x}}}_n^\star \cdot e^{-i 2 \pi \frac{k}{N_t} n}, 
\end{aligned}
\end{equation}
with $k = 0, \dots, N_t - 1$ and where 
$[{{\epsilon}_0^\star, {\epsilon}_1^\star,......{\epsilon}_{N_t-1}^\star}]$ and $[{\widetilde{\bm{x}}_0^\star, \widetilde{\bm{x}}_1^\star,......\widetilde{\bm{x}}_{N_t-1}^\star}]$ are the discrete Fourier coefficients.
Then the Power Spectral Density (PSD) is given by 

\begin{equation}\label{eq:PSD}
\begin{aligned}
\text {PSD}_{\epsilon} (k) = \dfrac{1}{N_t} |\hat{\epsilon}_{k}|^2, \quad \quad \text {PSD}_{\widetilde{\bm{x}}} (k) = \dfrac{1}{N_t} |{\hat{\bm{x}}_{k}}|^2.
\end{aligned}
\end{equation}

A PSD threshold is established to remove from the data the information associated to the frequencies that fall below such a value. 
It is worth noting that the PSD threshold needs to be properly set in order to keep within the data the frequencies containing the essential information associated to the underlying physical phenomena of the system. 
Its optimal value could depend on the problem at hand and there is no general rule to choose it. A high value may remove some significant frequencies leading to
a loss of important information. Conversely, a low value may not filter properly the snapshots by offering no substantial gain as it does not reduce the modal content of the system. 
In this work, the PSD
threshold value is determined through a trial and error approach. 
Of course, such a procedure is time-consuming and leaves room for further research for
improvement. 

\begin{figure}[ht] 
\centering
{\label{fig:Flowchart}\includegraphics[width=1\linewidth]{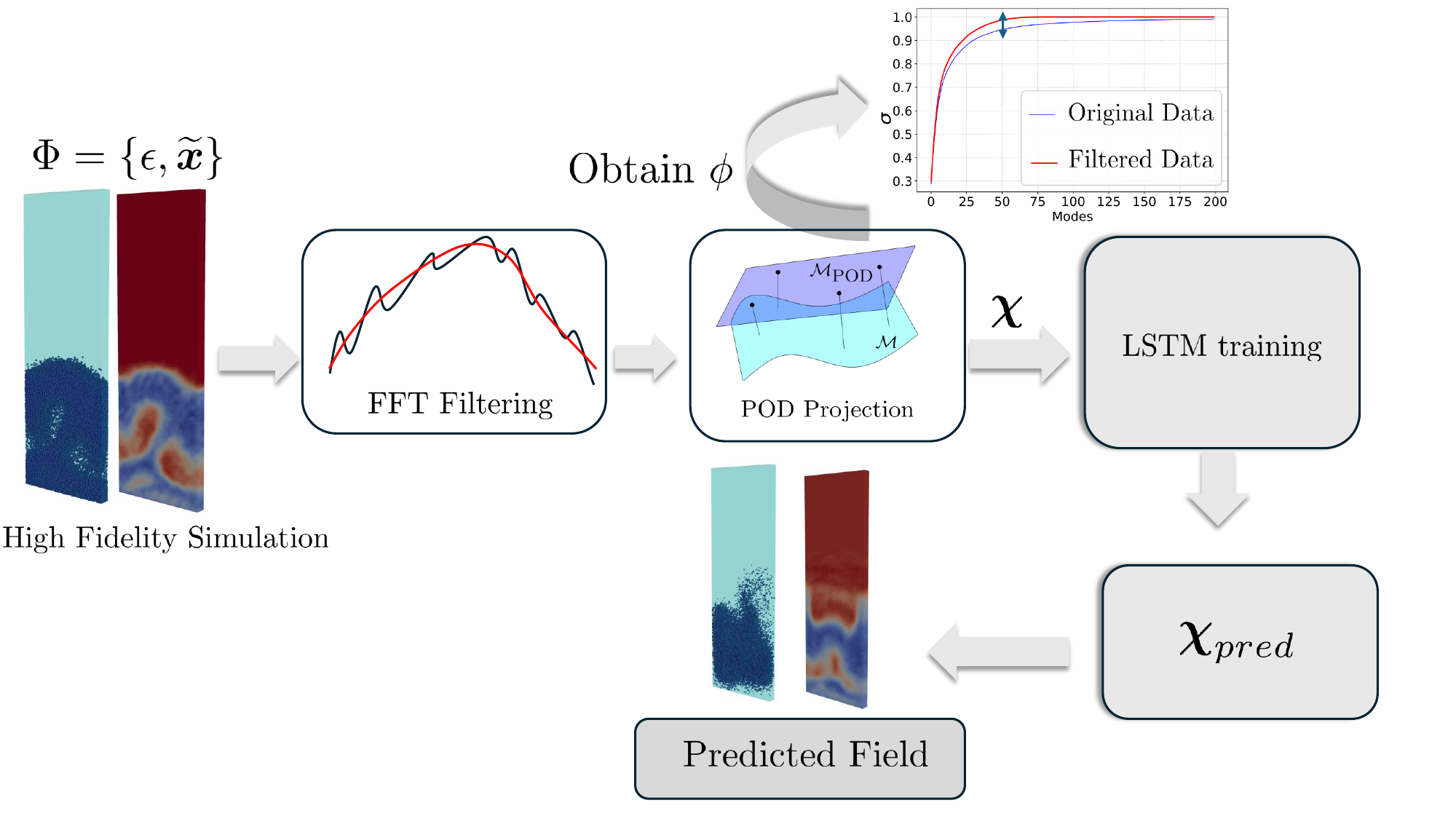}}
\caption{ 
Flowchart of the proposed ROM algorithm. 
} 
\label{fig:flowchart}
\end{figure}

Finally, Figure \ref{fig:flowchart} shows a sketch of the proposed ROM approach including the main steps of the offline and online stages. 



\section{Numerical results}\label{sec:res}
\graphicspath{{./img/}}

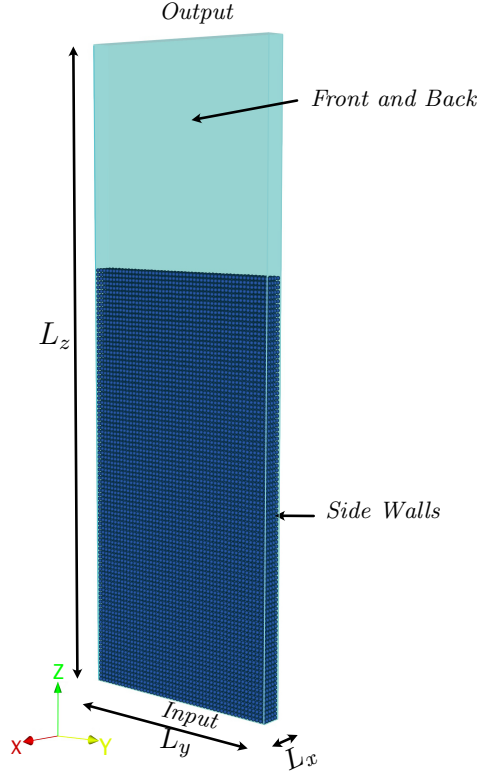
\begin{figure}
\centering
\begin{overpic}[width=0.3\linewidth, grid=false]{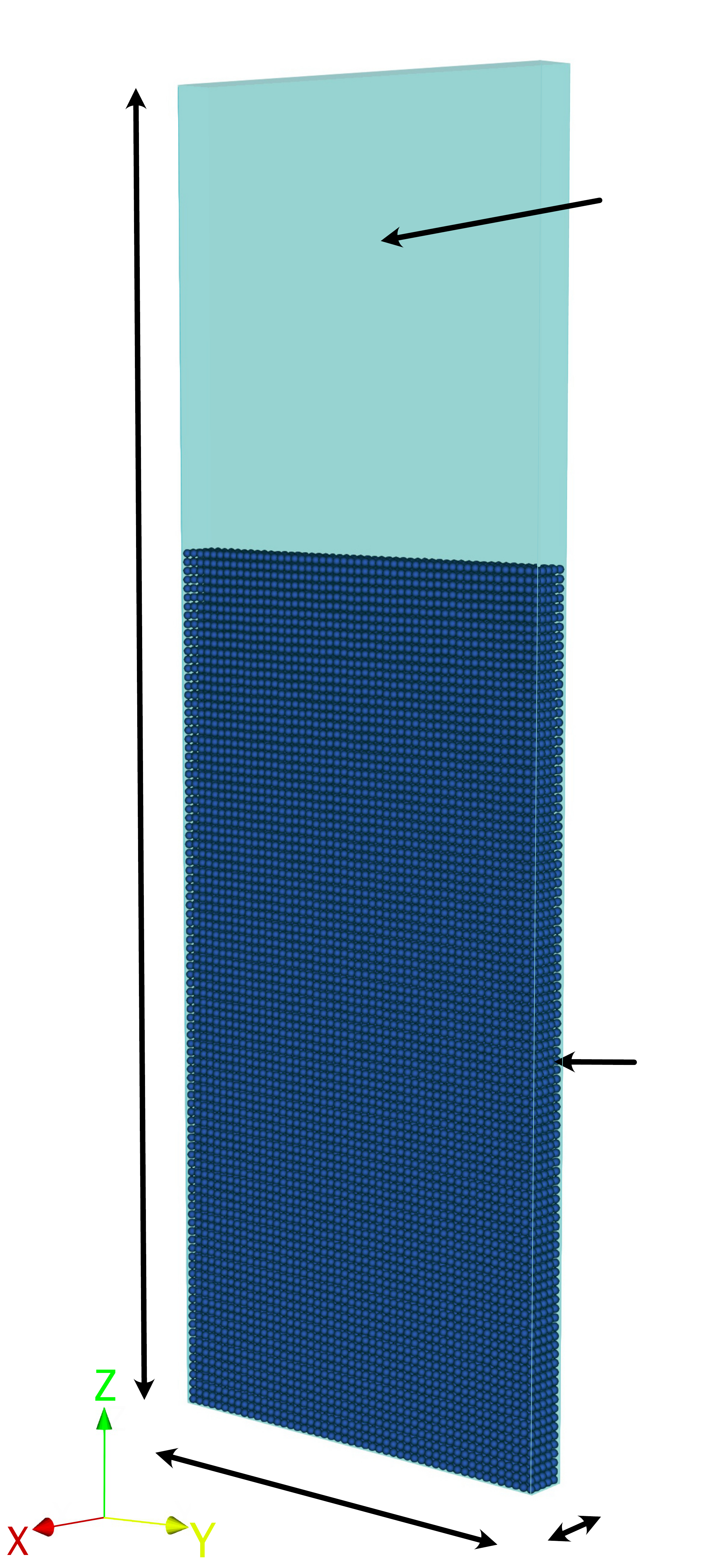}  
    \put(4,55){\Large{$L_z$}}
    \put(20,1.5){\large{$L_y$}}
    \put(36,-1.5){\rotatebox{25}{\large{$L_x$}}}
    \put(20,5.8){\rotatebox{-13}{\small{\emph{Input}}}}
    \put(20,98){\rotatebox{0}{\small{\emph{Output}}}}
    \put(42,32){\small{\emph{Side Walls}}}
    \put(40,86.5){\small{\emph{Front and Back}}}
\end{overpic}
    \caption{Sketch of the computational domain at $t = t_0 = 0$. The system is at rest at the beginning of the simulation. }
 \label{fig:BC}
\end{figure}

The performance of our proposed ROM approach have been evaluated against a gas-solid fluidized bed benchmark. In this two-phase flow system, the particles are transported by a carrier gas flow \cite{goldschmidt2001hydrodynamic, fernandes}. The computational domain for our simulation consists of a rectangle with dimensions $L_x \times L_y \times L_z = 15 \times 150 \times 450$ mm as depicted in \fig{fig:BC}. It is discretized into $N_x \times N_y \times N_z = 2 \times 30 \times 90$ cells. While the   motion equations of the particles are solved in 3D, the gas dynamics is simulated in a 2D framework in the $y-z$ plane \cite{goldschmidt2001hydrodynamic}.

Concerning the boundary
conditions, we refer to the labels reported in \fig{fig:BC} indicating the boundaries of our computational domain. 
A no-slip boundary condition is applied to the \emph{Side Walls} faces while the \emph{Front And Back} faces in the $y$-$z$ plane are set to symmetry. A mixed condition is employed to the \emph{Outlet} boundary: this adopts a homogeneous Neumann boundary condition in case of outflow, otherwise a null normal velocity to prevent backflow. Finally, a non-homogeneous Dirichlet boundary condition is applied on the \emph{Inlet} face: here, the interstitial velocity, which is the upward fluid flow velocity through the gaps between particles, is computed by dividing a pre-defined velocity (1.875 m/s) by the volume fraction.

A total number of $n_p = 24750$ spherical particles with $d_p = 0.0025$ m and $\rho_p = 2488.32$ m$^3$ are distributed uniformly in the domain.  The initial bed height is set to $300$ mm. 
The viscosity and density of the fluid phase are $\mu_f = 1e-5$ Pa $\cdot$ s and $\rho_f = 1000$ Kg/m$^3$, respectively. This results in a Stokes number $Stk = 300$ (see eq. \eqref{Stk}).
The system is initially at rest. The simulation is run up to $T = 5$ s by adopting a time-step $\Delta t = 2e-5$ s keeping the maximum Courant number below one for the sake of numerical stability.

\begin{figure}\vspace{2cm}
    \centering
    
\begin{subfigure}{.6\textwidth}
  \begin{overpic}[width=\linewidth, grid=false]{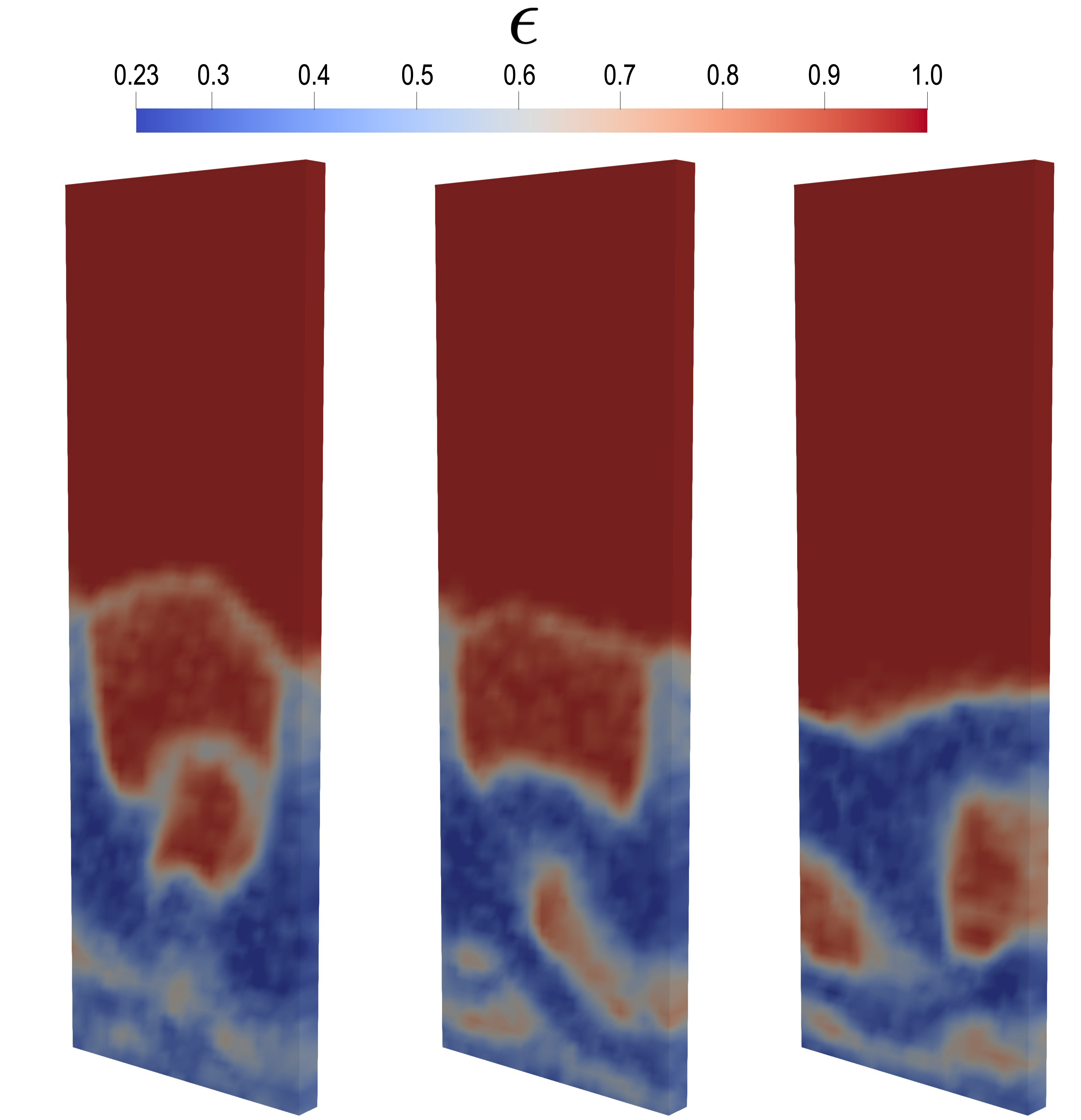}  
    \put(14,80){\tiny{\textcolor{white}{$t$=2 s}}}
    \put(46,80){\tiny{\textcolor{white}{$t$=3.5 s}}}
    \put(80,80){\tiny{\textcolor{white}{$t$=5 s}}}
  \end{overpic}
  \caption{Time evolution of the fluid volume fraction: $t = 2$ s (left), $t = 3.5$ s (center) and $t = 5$ s (right).}
 \label{fig:Eurlerian_FOMSol}
\end{subfigure}

\medskip

\begin{subfigure}{.6\textwidth}
  \begin{overpic}[width=\linewidth, grid=false]{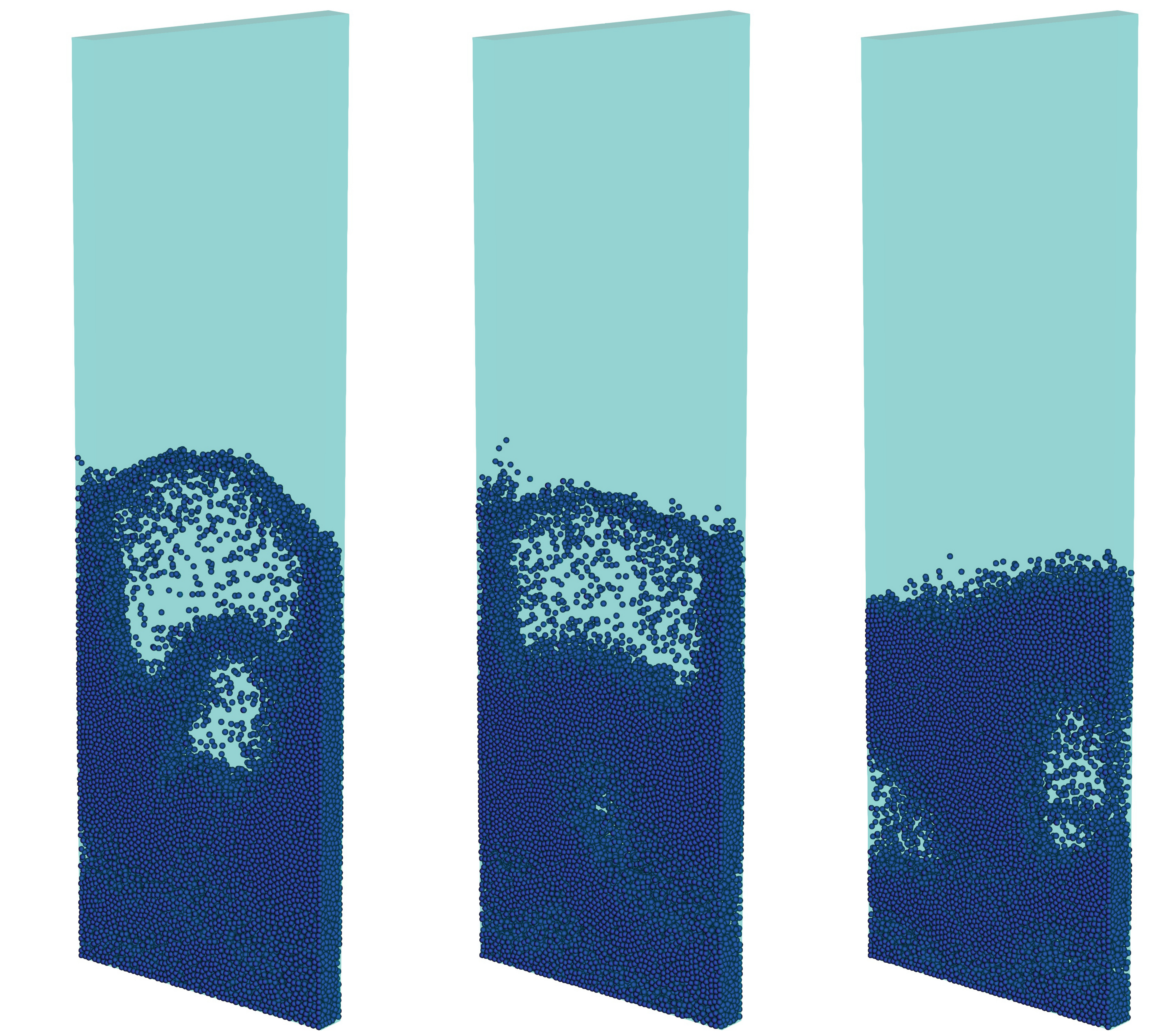}  
    \put(14,82){\tiny{\textcolor{white}{$t$=2 s}}}
    \put(47,82){\tiny{\textcolor{white}{$t$=3.5 s}}}
    \put(82,82){\tiny{\textcolor{white}{$t$=5 s}}}
  \end{overpic}
  \caption{Time evolution of the particles distribution: $t = 2$ s (left), $t = 3.5$ s (center) and $t = 5$ s (right).}
  \label{fig:Lagrang_FOMSol}
\end{subfigure}

    \caption{FOM solution of the fluid volume fraction (A) and particle position (B)  at times $t = 2$ s (initial  column), $t = 3.5$ s (middle column) and $t = 5$ s (last column).}
\label{fig:full_order}
\end{figure}
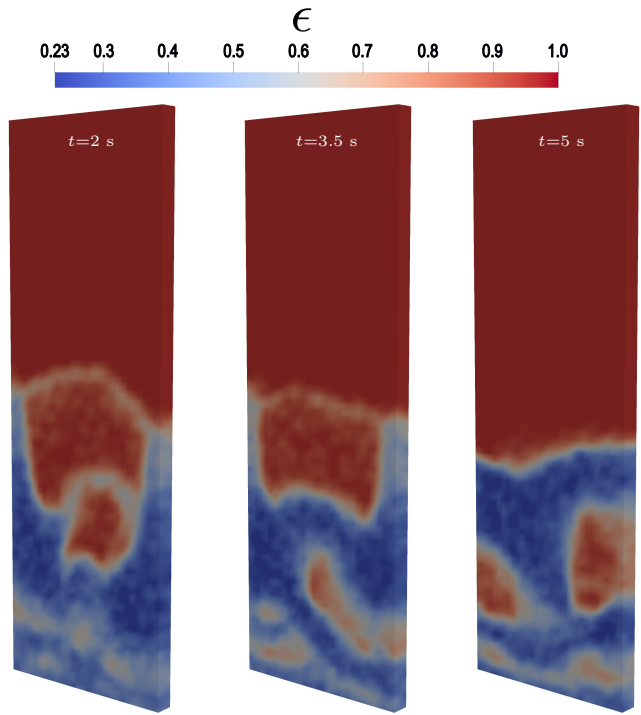
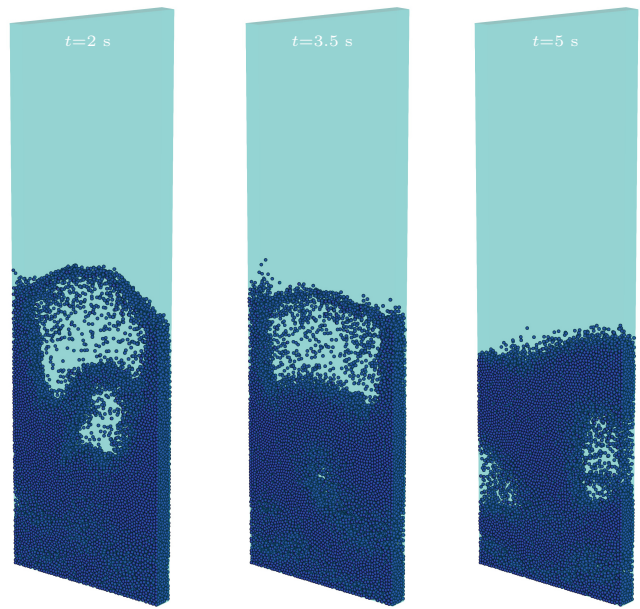

\fig{fig:full_order} demonstrates the time evolution of the fluid-particle system over (0, 5] s: full order solutions for the void fraction and particle position are depicted at $t = 2$ s, $t = 3.5$ s and $t = 5$ s from left to right. As expected, the particles start their motion by the upward gas flow and suspend, subsequently, when the gas reaches the fluidization velocity \cite{goldschmidt2001hydrodynamic, fernandes, hajisharifi2023non}.  

Firstly we are going to investigate the performance in terms of accuracy of the ROM model in the identification and prediction of the behaviour of the flow field in Sec. \ref{sec:Eulerian_rec}. Then the results related to particle position are discussed in Sec. \ref{sec:Lag_rec}. Finally some insights about the efficiency of our approach are given in Sec. \ref{sec:cost}.

\subsection{Eulerian phase}\label{sec:Eulerian_rec}
We collect an original database over (0, 5] s consisting of 500 snapshots with a sampling frequency of 100 Hz, i.e. every 0.01 seconds.
These snapshots are divided into two different sets. A first set, called training set, is used to generate
the reduced basis. 
The second set, called validation set, is the complement of the training set in the original
database and it is used to assess the accuracy of the ROM solution. 
Out of the 500 computed $\epsilon$ in the original database, we take 450 (i.e., 90\% of the database)
to form the training set. These 450 solutions are the first
450 in the database (associated to the time interval (0, 4.5]) . 
The remaining 50 solutions form the validation set. Concerning the LSTM network, it consists of one only layer and 1024 neurons, the learning rate is
$2e-4$, the number of epochs is 3000 corresponding to a MSE final value of $4e-4$. Moreover, we set $s = 50$.  

\begin{figure}[t]\centering
\subfloat[PSD versus frequency for the grid point $P_1$]{\label{fig:50PercentTraining}\includegraphics[width=.53\linewidth]{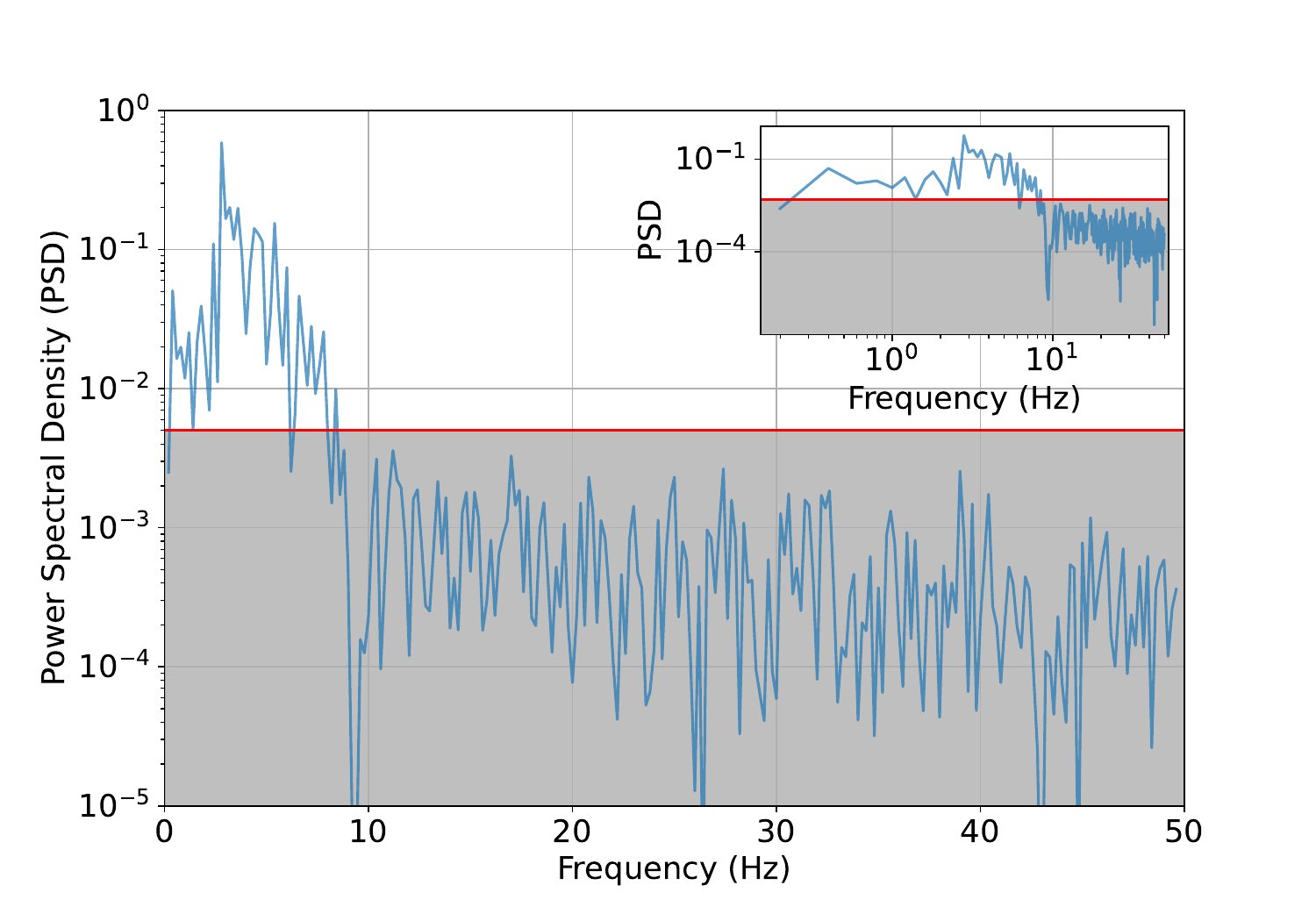}}
\subfloat[PSD versus frequency for the grid point $P_2$]{\label{fig:70PercentTraining}\includegraphics[width=.53\linewidth]{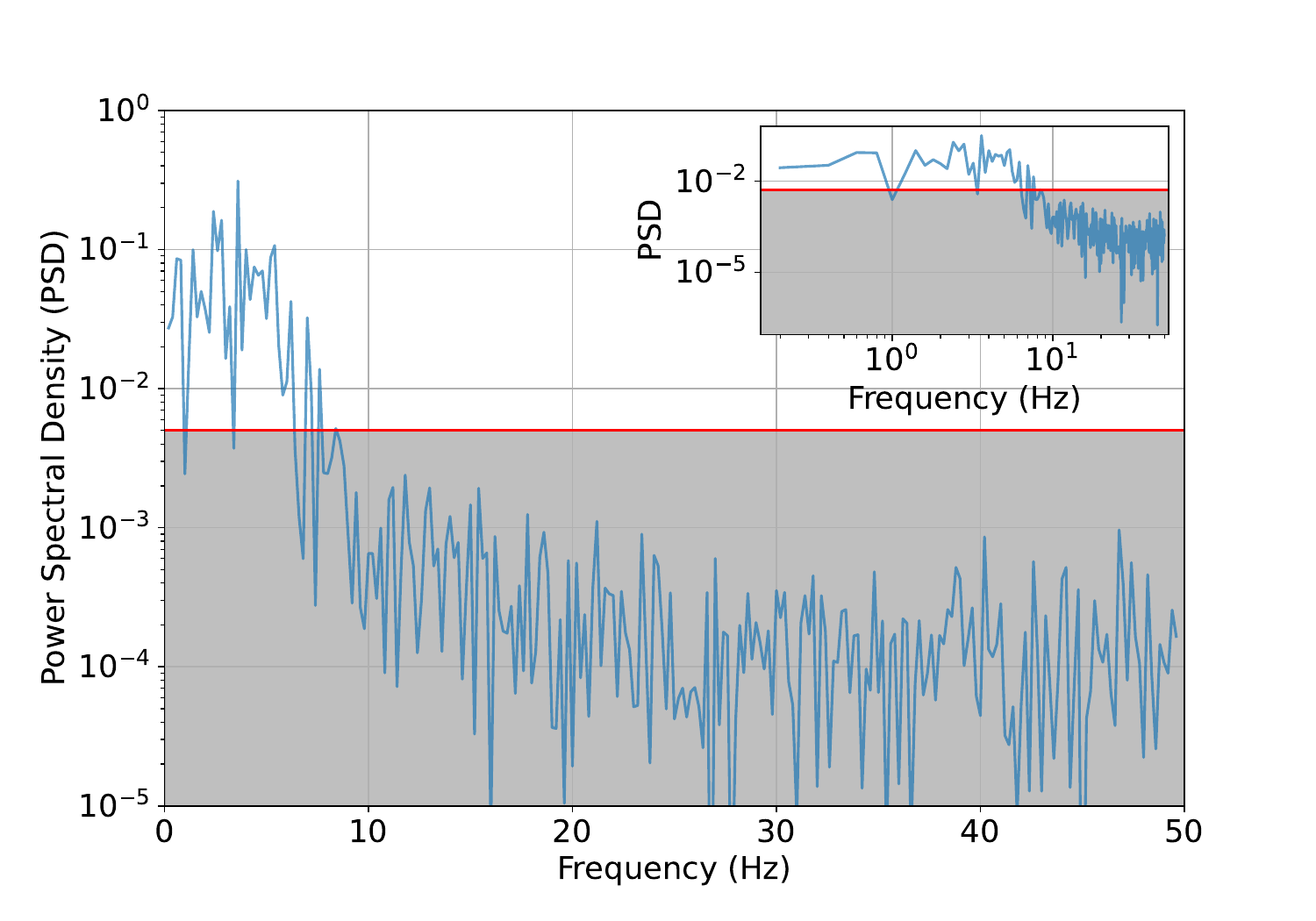}}\par 
\subfloat[PSD versus frequency for the grid point $P_3$]{\label{fig:90PercentTraining}\includegraphics[width=.53\linewidth]{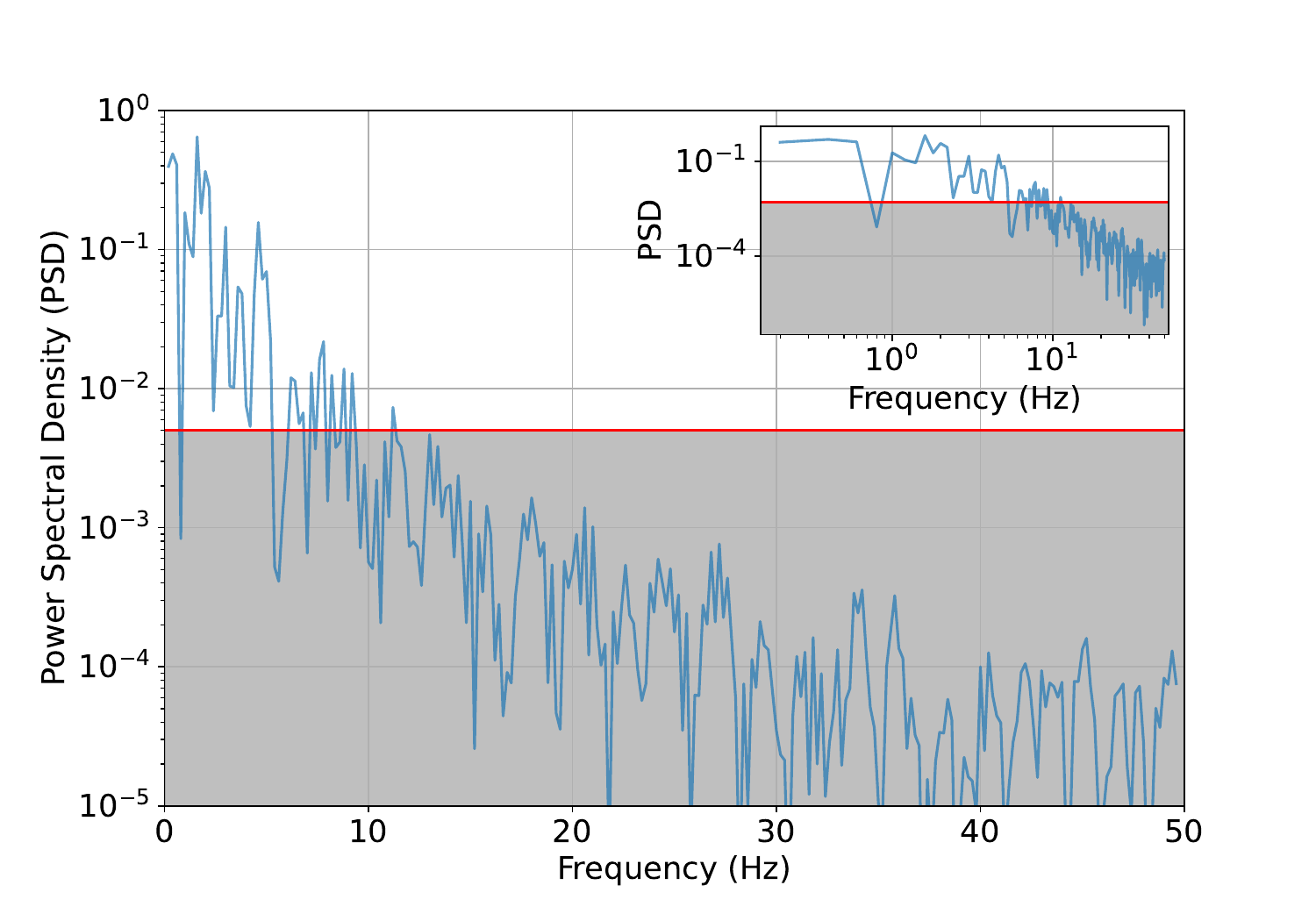}}
\subfloat[PSD versus frequency for the grid point $P_4$]{\label{fig:90PercentTraining}\includegraphics[width=.53\linewidth]{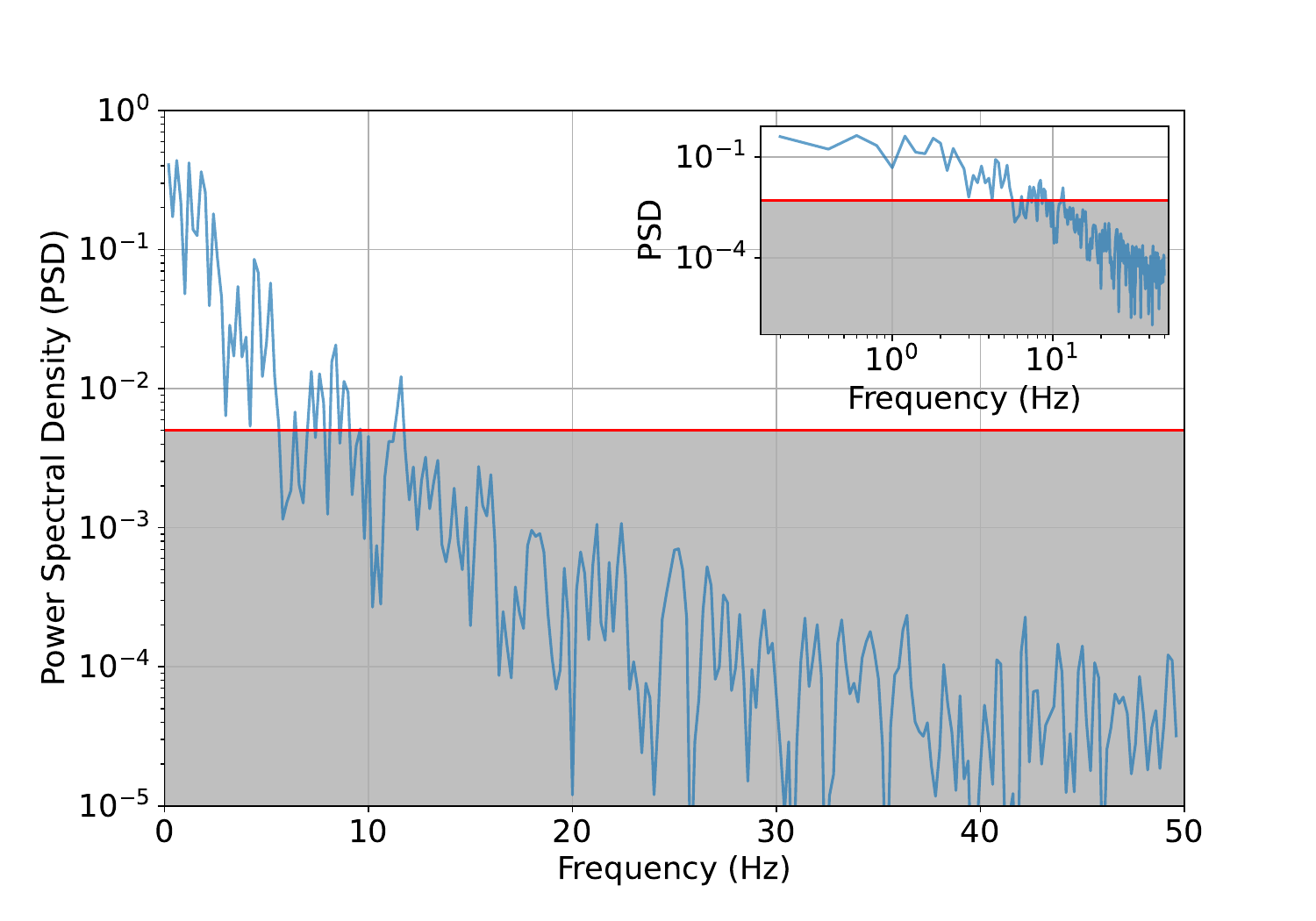}}

\caption{Filtering validation - Eulerian field: PSD analysis across four distinct grid points: $P_1 = ($0, -0.07, 0.025$)$  m, $P_2 =($0, 0.07, 0.025$)$ m, $P_3 = ($0, -0.07, 0.25$)$ m and $P_4 = ($0, -0.07, 0.25$)$ m. The main figure is in logarithmic scale along $y$ axis while the inset is in logarithmic scale for both $x$ and $y$ axes. The red line shows the PSD threshold (equal to 0.005) and the cut zone is highlighted in grey.} 
\label{fig:PSD}
\end{figure}

\begin{figure}
\vspace{1cm}

\begin{subfigure}{.65\textwidth}
  \begin{overpic}[width=\linewidth, grid=false]{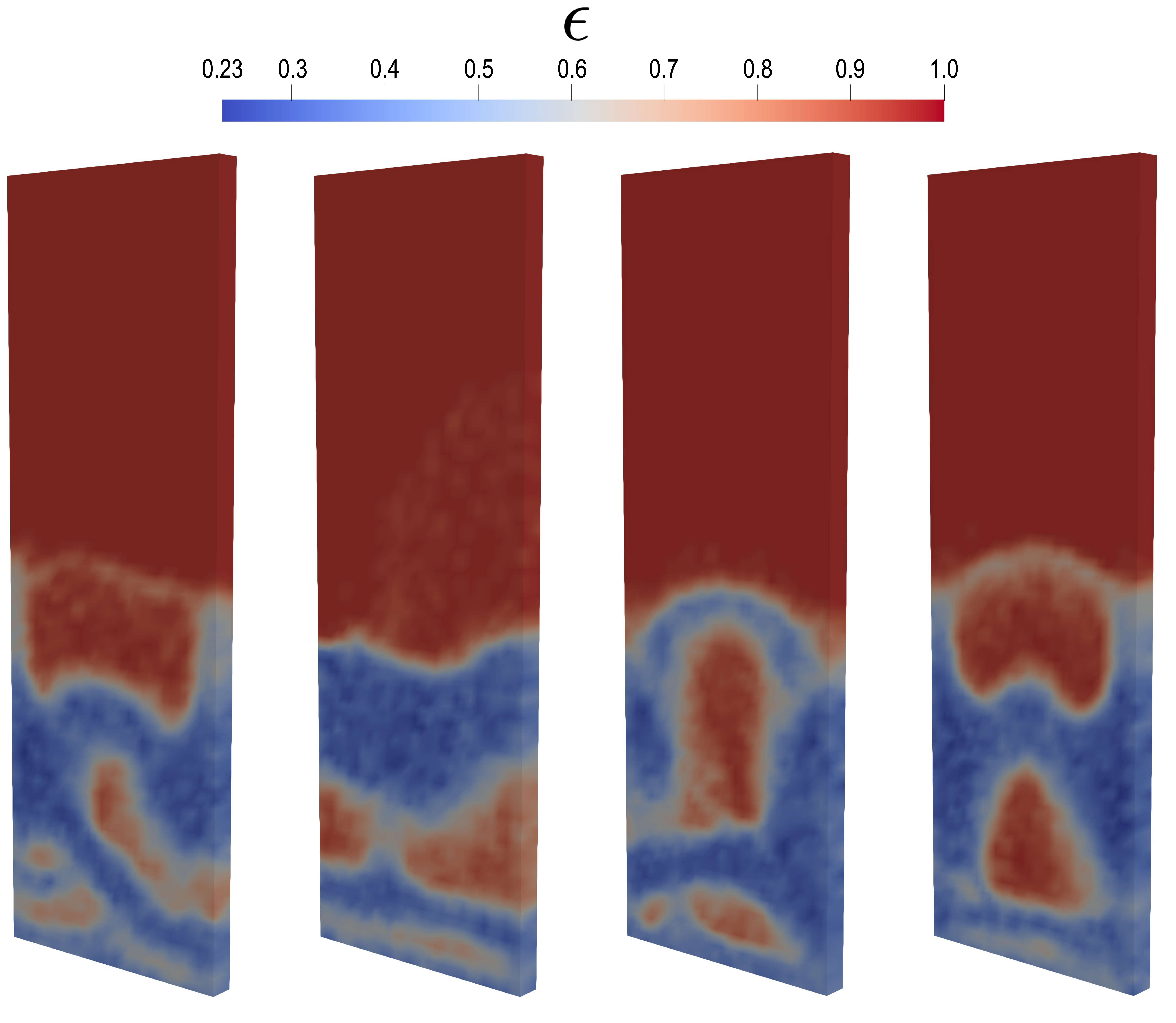}  
    \put(5,65){\tiny{\textcolor{white}{$t$=3.5 s}}}
    \put(32,65){\tiny{\textcolor{white}{$t$=4.6 s}}}
    \put(58,65){\tiny{\textcolor{white}{$t$=4.7 s}}}
    \put(84,65){\tiny{\textcolor{white}{$t$=4.8 s}}}
  \end{overpic}
  \caption{Unfiltered FOM solution at different time instances: $t=3.5$ s, $t=4.6$ s, $t=4.7$ s and $t=4.8$ s from left to right, respectively.}
\end{subfigure}
\par
\vspace{0.5cm}
\begin{subfigure}{.65\textwidth}
  \begin{overpic}[width=\linewidth, grid=false]{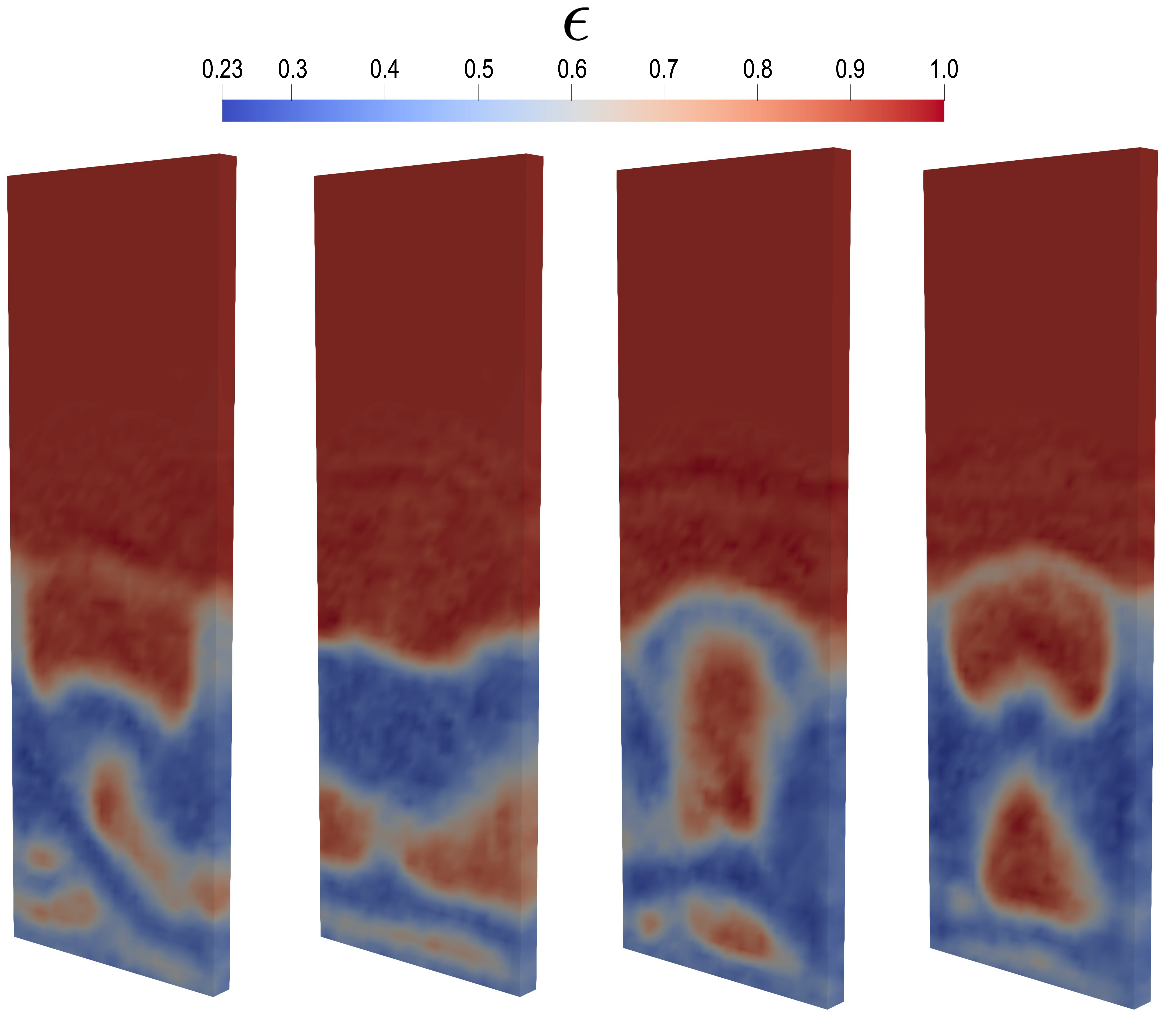}  
    \put(5,65){\tiny{\textcolor{white}{$t$=3.5 s}}}
    \put(32,65){\tiny{\textcolor{white}{$t$=4.6 s}}}
    \put(58,65){\tiny{\textcolor{white}{$t$=4.7 s}}}
    \put(84,65){\tiny{\textcolor{white}{$t$=4.8 s}}}
  \end{overpic}
  \caption{Filtered FOM solution at different time instances: $t=3.5$ s, $t=4.6$ s, $t=4.7$ s and $t=4.8$ s from left to right, respectively.}
\end{subfigure}

\caption{Filtering validation - Eulerian field: comparison of the time evolution of $\epsilon$ computed by unfiltered (first row) and filtered (second row) FOM.}
\label{fig:Vis_FOM_filtered_Eulerian}
\end{figure}

We start to illustrate the filtering process. \fig{fig:PSD} shows the PSD as a function of frequency for four different grid points: $P_1 = ($0, -0.07, 0.025$)$  m, $P_2 =($0, 0.07, 0.025$)$ m, $P_3 = ($0, -0.07, 0.25$)$ m and $P_4 = ($0, -0.07, 0.25$)$ m. For each point, the main figure is in  logarithmic scale along $y$-axis.  On the other hand, 
the inset is in  logarithmic scale for both $x$ and $y$ axes. 
The PSD threshold value, represented by the red horizontal line, is set to 0.005. 
We can see that the frequency content greater than about 10 Hz is cut. 
In order to validate such a setting, \fig{fig:Vis_FOM_filtered_Eulerian} shows a qualitative comparison between the original FOM snapshots and the filtered ones (i.e., the FOM snapshots resulting from the filtering process). We observe that unfiltered and filtered FOM solutions are in perfect agreeement confirming  that the selected PSD threshold is well-suited for the system at hand. 


\begin{figure} 
    \centering  \includegraphics[width=90mm,scale=0.5]{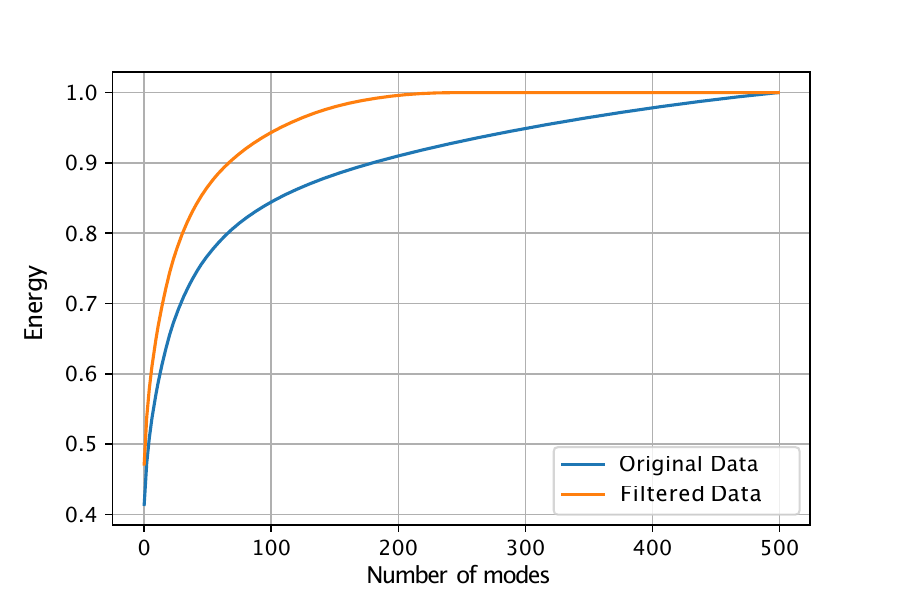}
        \caption{ROM validation - Eulerian field: comparison of the cumulative energy $E$ defined in eq. \eqref{k_equation} of the fluid volume fraction $\epsilon$ for the original and filtered FOM snapshots. }
\label{fig:Energy}
\end{figure}

For further validation we compare the POD analysis of the filtered and original data. 
We plot the cumulative energy $E$ (eq. \eqref{k_equation}) in \fig{fig:Energy}. 
It clearly highlights the impact of the filtering process: the cumulative energy associated to the unfiltered snapshots shows a significantly slower convergence with respect to the filtered one. 
In order to provide a more quantitative comparison, Tab.~\ref{tab:Modes_num} reports the number of modes associated to $\delta = 0.7, 0.9, 0.99$: 
we can see that by performing the POD on the original dataset, more than double the number of modes are required for $\delta = 90\%$ and $99\%$. 

\begin{table}[ht!]
\centering
\begin{tabular}{|c|c|c|}
\hline  
   & Original FOM data  & Filtered FOM data  \\ \hline
$\delta = 70 \%$ & 29  &  15   \\ \hline
$\delta = 90 \%$ & 180  & 67  \\ \hline
$\delta = 99 \%$ & 449 &  176  \\\hline 
\end{tabular}
  \caption{ROM validation - Eulerian field: number of required modes to retain different energy thresholds, $\delta = 0.7, 0.9, 0.99$, for original and filtered FOM data.
  }
 \label{tab:Modes_num}
\end{table}

In our study we set $\delta = 90 \%$ corresponding to 67 POD modes. This choice is a trade-off that ensures a good balance between accuracy and efficiency. 
Infact, while a higher energy threshold could potentially increase the accuracy of the ROM solution, it badly affects the performance of LSTM which is less efficient in predicting higher dimensional data. 

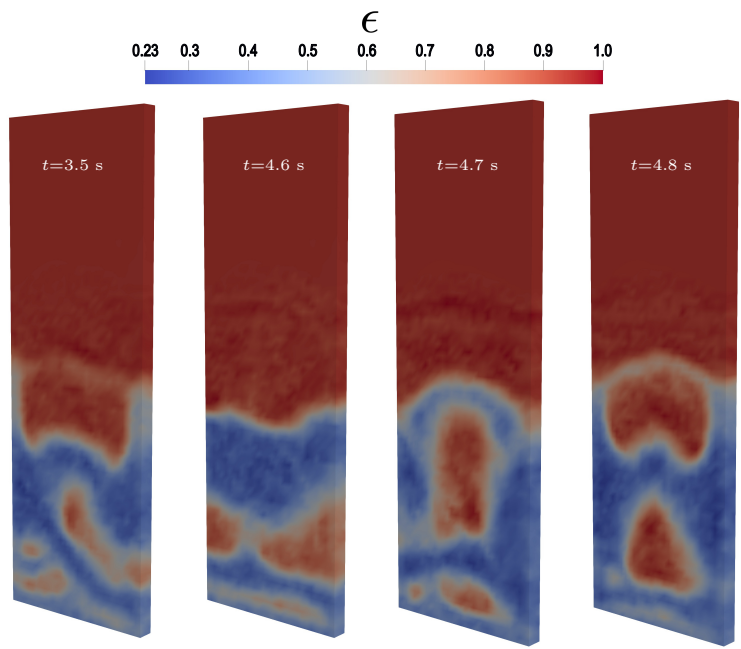
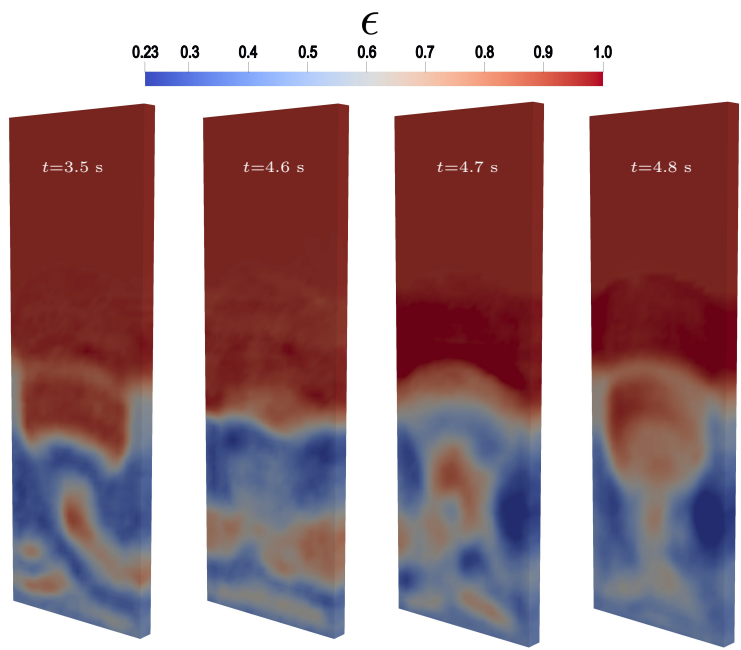
\begin{figure}
\vspace{1cm}

\begin{subfigure}{.65\textwidth}
  \begin{overpic}[width=\linewidth, grid=false]{total_FOM_Filtered.pdf}  
    \put(5,65){\tiny{\textcolor{white}{$t$=3.5 s}}}
    \put(32,65){\tiny{\textcolor{white}{$t$=4.6 s}}}
    \put(58,65){\tiny{\textcolor{white}{$t$=4.7 s}}}
    \put(84,65){\tiny{\textcolor{white}{$t$=4.8 s}}}
  \end{overpic}
  \caption{FOM solution at different time instances: $t=3.5$ s, $t=4.6$ s, $t=4.7$ s and $t=4.8$ s from left to right, respectively.}
\end{subfigure}
\par
\vspace{0.5cm}
\begin{subfigure}{.65\textwidth}
  \begin{overpic}[width=\linewidth, grid=false]{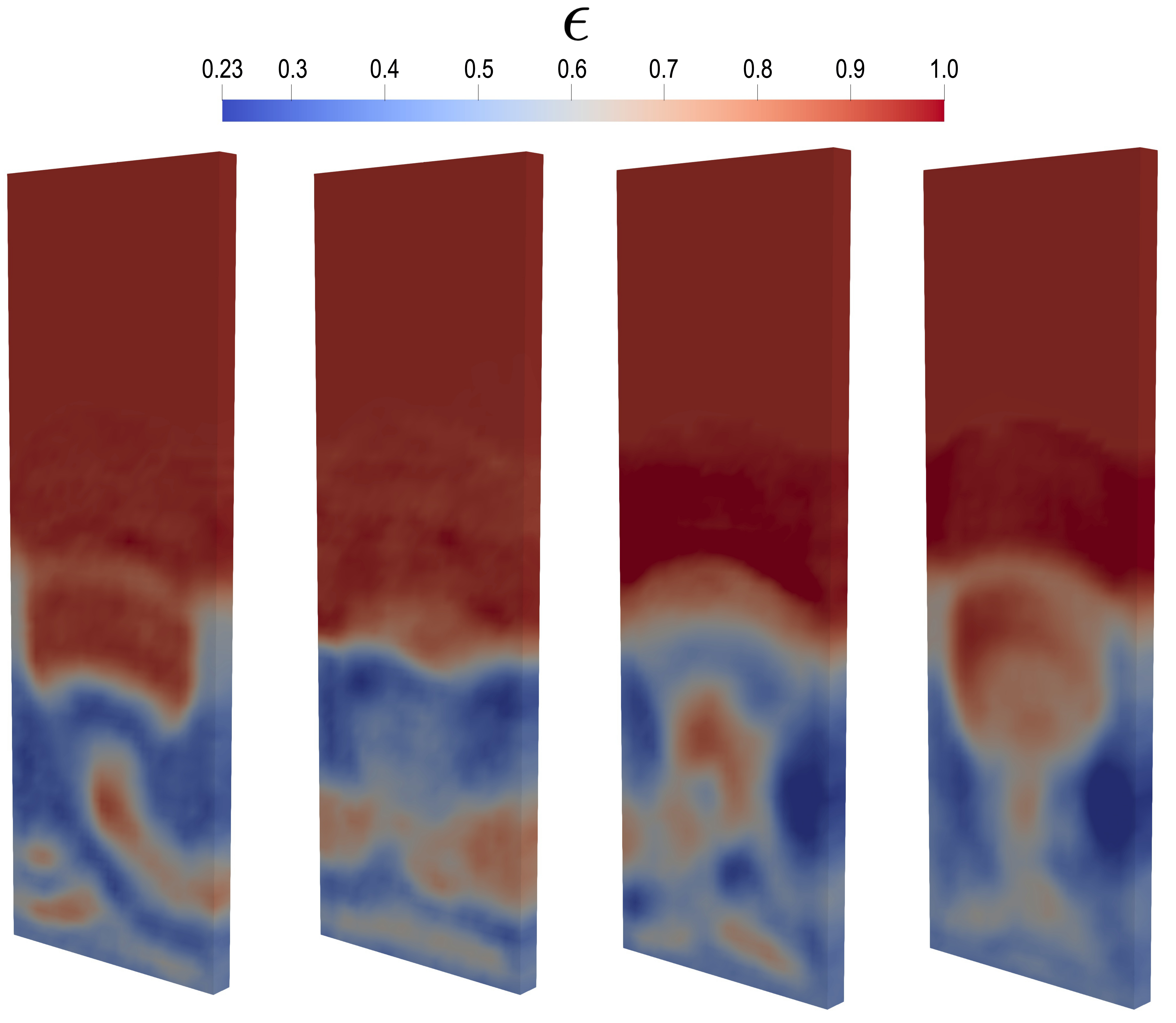}  
    \put(5,65){\tiny{\textcolor{white}{$t$=3.5 s}}}
    \put(32,65){\tiny{\textcolor{white}{$t$=4.6 s}}}
    \put(58,65){\tiny{\textcolor{white}{$t$=4.7 s}}}
    \put(84,65){\tiny{\textcolor{white}{$t$=4.8 s}}}
  \end{overpic}
  \caption{ROM solution at different time instances: $t=3.5$ s, $t=4.6$ s, $t=4.7$ s and $t=4.8$ s from left to right, respectively.}
\end{subfigure}

\caption{ROM validation - Eulerian field: comparison of the time evolution of $\epsilon$ computed by FOM (first row)  and by ROM (second row) for $\delta = 90\%$. Notice that $t = 3.5$ s belongs to the training set while the other time instances are in the validation set.}
\label{fig:Vis_Filtered_ROM_Eulerian}
\end{figure}

A qualitative comparison of (filtered) FOM and ROM solutions is depicted in \fig{fig:Vis_Filtered_ROM_Eulerian} for $t = 3.5, \ 4.5, \ 4.7$ and 4.8 s. Note that $t = 3.5$ s belongs to the training set, so it is depicted to evaluate the ROM capability to identify the system dynamics. The remaining
three times, i.e., $t = 4.5, \ 4.7$ and $4.8$ s, are not associated with the training set and thus are used to check the accuracy of the ROM in predicting the system dynamics. Let us discuss the results in \fig{fig:Vis_Filtered_ROM_Eulerian} starting from the system identification. We observe that our ROM approach is able to provide a very accurate reconstruction of the solution for $t=3.5$ s. 
Now, let us take a look at the solutions corresponding to the time not associated with the training set. 
The visual comparison for $t=4.6$ s 
demonstrates no significant difference between the FOM and ROM predicted solution. 
However,  for $t=4.7$ and $t=4.8$, although ROM is of course able to forecast the main features of the dynamics, such as the bed height, it cannot capture all the flow structures.  
This comparison suggests a high accuracy level for system identification and short-time prediction. 
On the other hand, for longer time prediction, the model shows limitations in predicting the detailed dynamics of the system although from an engineering standpoint the results are still acceptable. 
\begin{figure} 
    \centering  \includegraphics[width=80mm,scale=0.5]{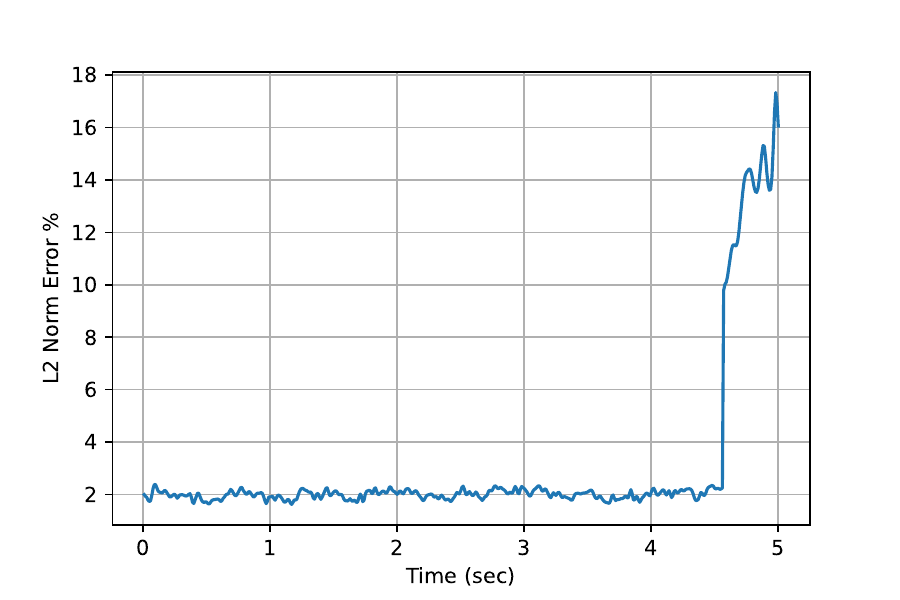}
        \caption{ROM validation - Eulerian field: time evolution of the $L^2$-norm relative error defined in eq. \eqref{eq:l2Error} for $\delta= 90 \%$.}
\label{fig:L2Error}
\end{figure}

In order to provide more quantitative results, we compute the time evolution of the $L^2$-norm error between ROM and FOM:

\begin{equation}
E_{\epsilon}(t) = 100 \cdot \dfrac{||\epsilon_h(t) - \epsilon_r(t)||_{L^2(\Omega)}}{||{\epsilon_h}(t)||_{L^2(\Omega)}},
\label{eq:l2Error}
\end{equation}
where $\epsilon_h$ is the fluid volume fraction computed with the FOM and $\epsilon_r$ is the corresponding field computed with the ROM.
In Fig. \ref{fig:L2Error} it is observed that the error remains around $2\%$ for the entire training set, confirming that ROM is accurate in performing system identification. This is aligned with the qualitative results shown in the first column of \fig{fig:Vis_Filtered_ROM_Eulerian}. On the other hand, we see that the error increases sharply around $t$ = 4.5 s up to $10\%$.  The magnitude of the error at the end of the time interval is around 17\%, confirming that while ROM is capable of short-time forecasting, its performance diminishes over a longer time period. However, these results are significantly better than the ones provided in \cite{li2023physics} where the authors adopted a Physics informed-DMD strategy. In \cite{li2023physics} the training set consists of 5001 snapshots collected in the interval [5, 10] s whilst 30 snapshots belonging to the interval [10.001, 10.030] s (which is almost 17 times smaller than ours) form the validation set: the relative error starts from about 7\% for $t = 10.001$, increases over the time and abundantly overcomes 30\% for $t = 10.030$ s.





\subsection{Lagrangian Phase}\label{sec:Lag_rec}

We collect an original database over [0, 4.6] s consisting of 460 snapshots with
a sampling frequency of 100 Hz, i.e. every 0.01 seconds, which are divided into a training set to generate the reduced basis and into a validation set to be used to assess the accuracy of the ROM solution. We take the first 450 snapshots (i.e., almost 98\% of the database) to form the training set (associated to the time interval (0, 4.5] s ) whilst the remaining 10 solutions form the validation set. It should be noted that since the Lagrangian dynamics is significantly more complex than the Eulerian one, we consider a smaller prediction window: [4.5, 4.6] s. Concerning the LSTM network, it consists of one only layer and $10$ neurons, the learning rate is 0.0001, the number of epochs is 3000 corresponding to a MSE final value of $1e-5$. Moreover, we set $s = 10$.  

\begin{figure}[h!]\centering
\subfloat[PSD versus frequency for $\widetilde{z}$ and particle \\ $l$ = 1000]{\label{fig:z1000}\includegraphics[width=.49\linewidth]{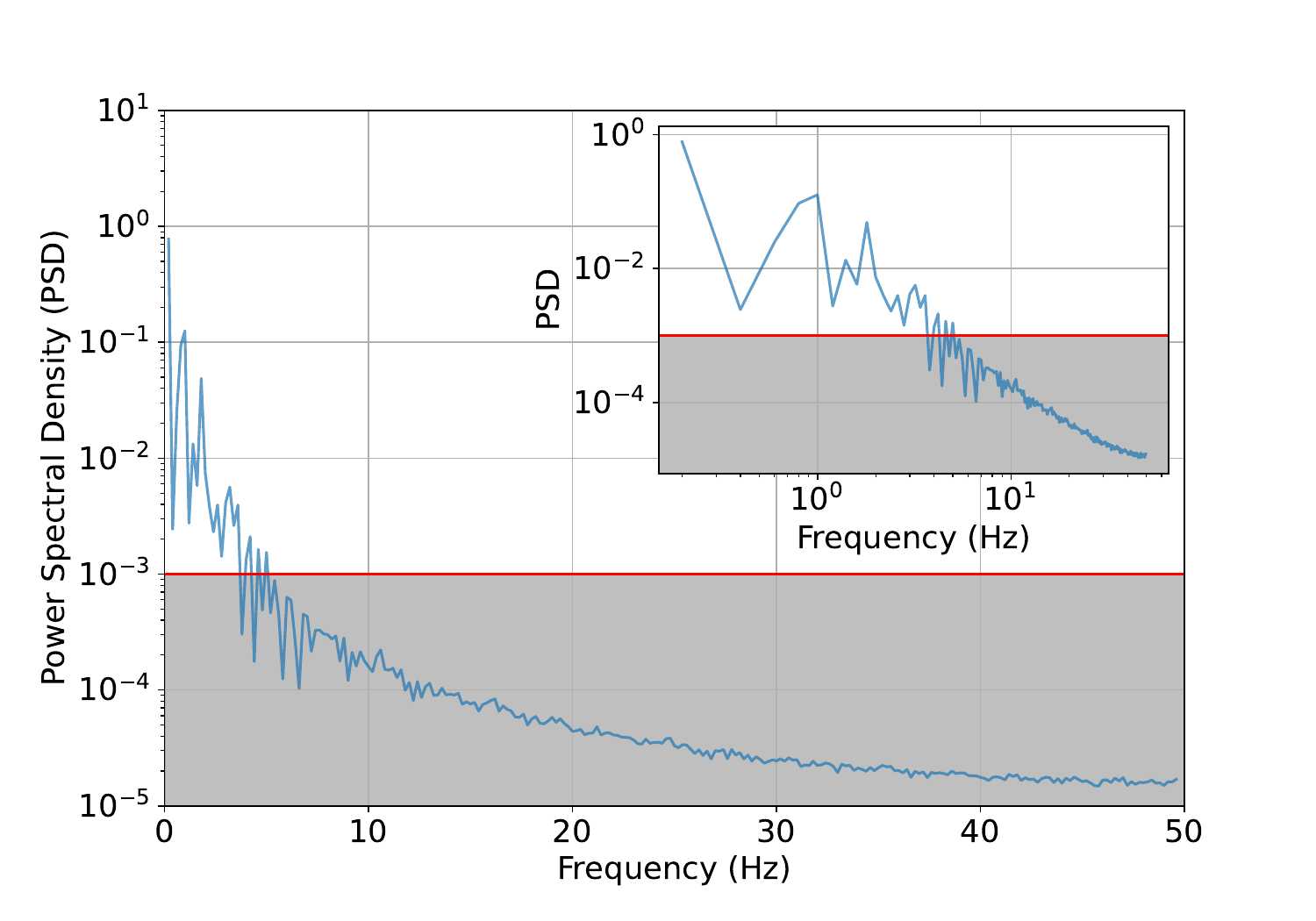}}
\subfloat[PSD versus frequency for $\widetilde{z}$ and particle \\ $l$ = 24000]{\label{fig:z24000}\includegraphics[width=.49\linewidth]{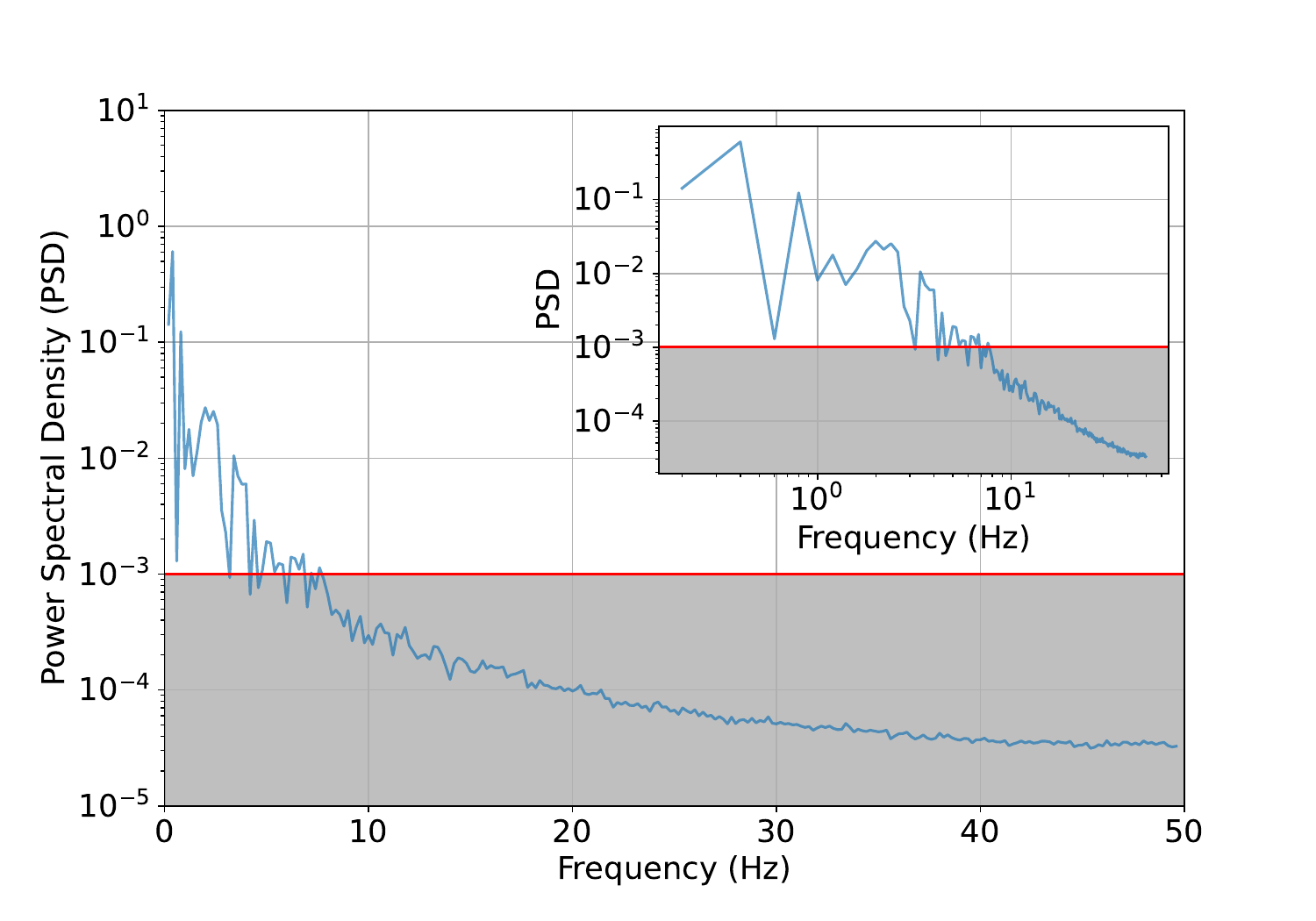}}\par 

\subfloat[PSD versus frequency for  $\widetilde{y}$ and particle \\ $l$ = 1000]{\label{fig:y1000}\includegraphics[width=.49\linewidth]{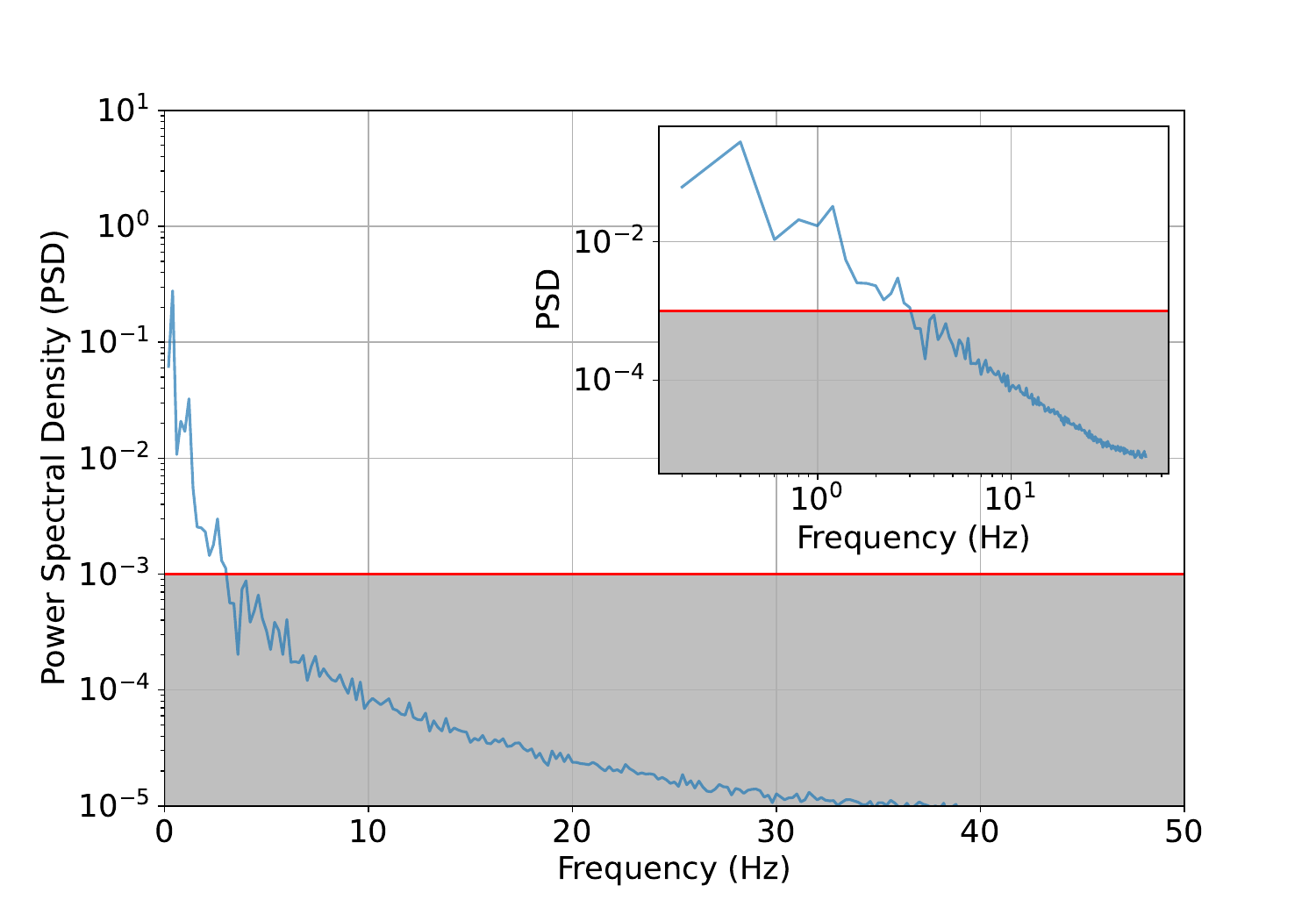}}
\subfloat[PSD versus frequency for  $\widetilde{y}$ and particle \\ $l$ = 24000]{\label{fig:y24000}\includegraphics[width=.49\linewidth]{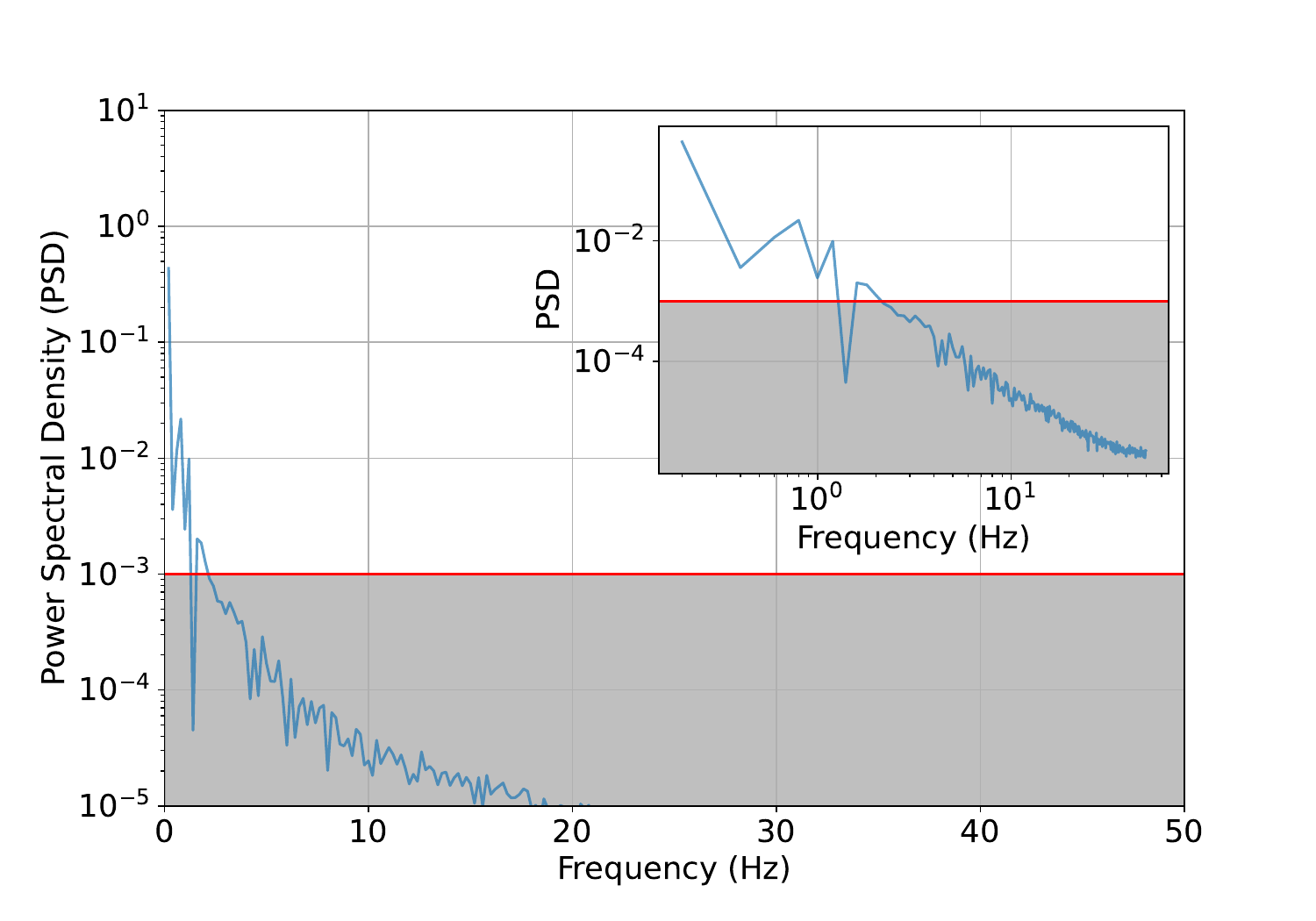}}\par

\subfloat[PSD versus frequency for $\widetilde{x}$ and particle \\ $l$ = 1000]{\label{fig:x1000}\includegraphics[width=.49\linewidth]{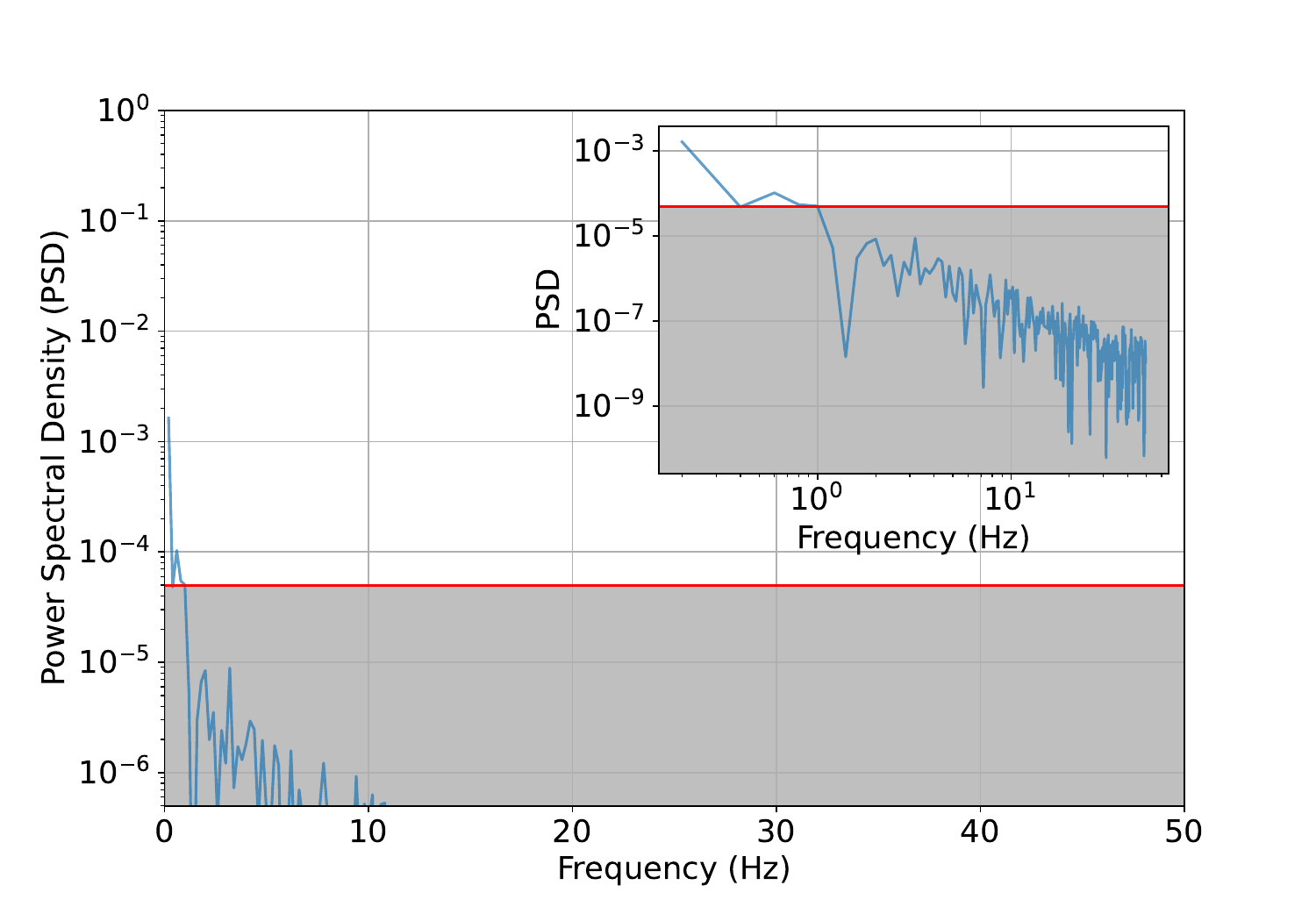}}
\subfloat[PSD versus frequency for $\widetilde{x}$ and particle \\ $l$ = 24000]{\label{fig:x24000}\includegraphics[width=.49\linewidth]{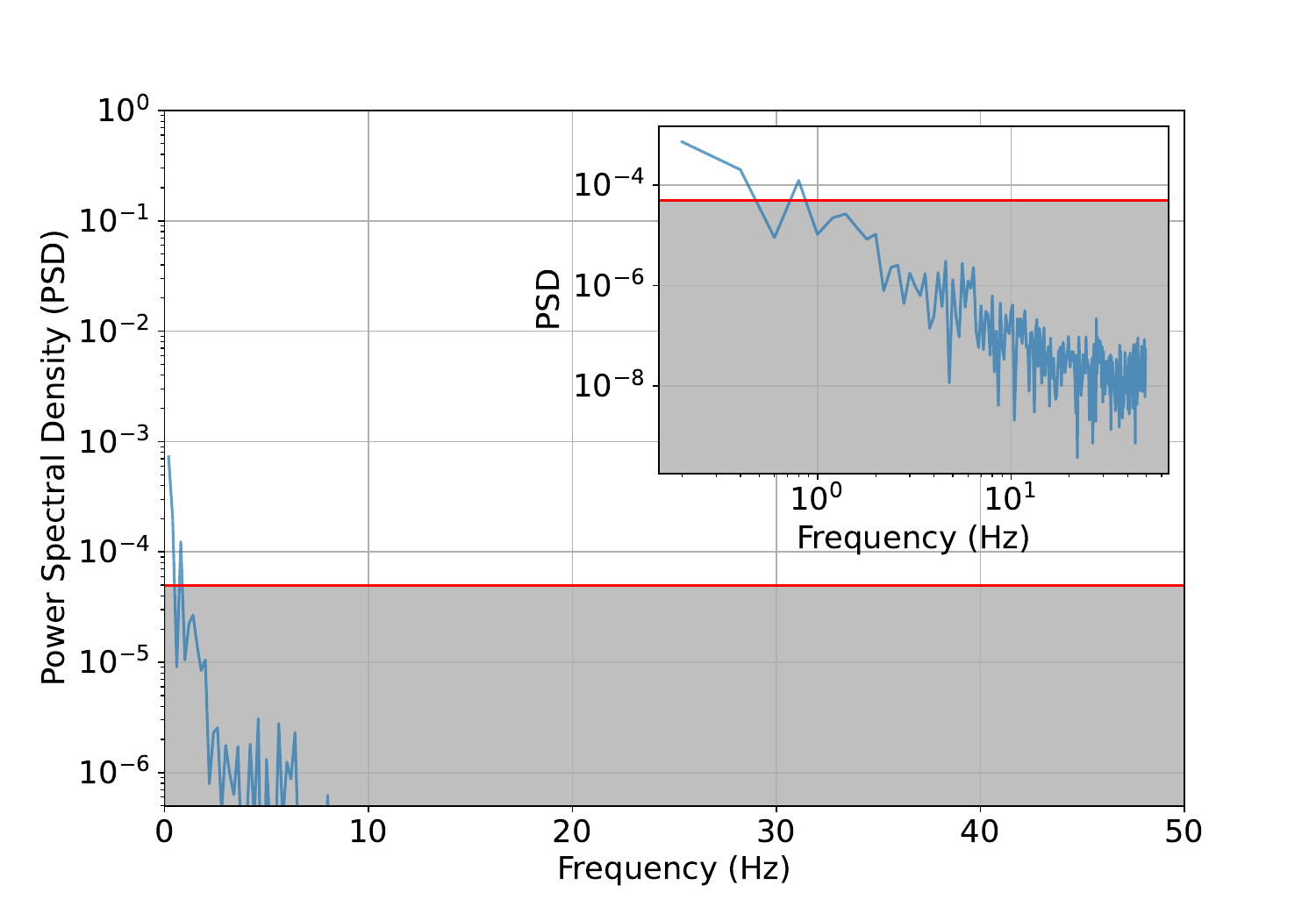}}

\caption{Filtering validation - Lagrangian field: PSD analysis across two distinct particles, $l = 1000$ and $l = 24000$. The main figure is in logarithmic scale along  $y$ axis while the inset is in logarithmic scale for both $x$ and $y$ axes. The red line shows the PSD threshold (equal to 0.001 for $\widetilde{y}$ and $\widetilde{z}$ and $5e^{-5}$ for $\widetilde{x}$) and the cut zone is highlighted in grey.} 
\label{fig:PSD_lagrangian}
\end{figure}

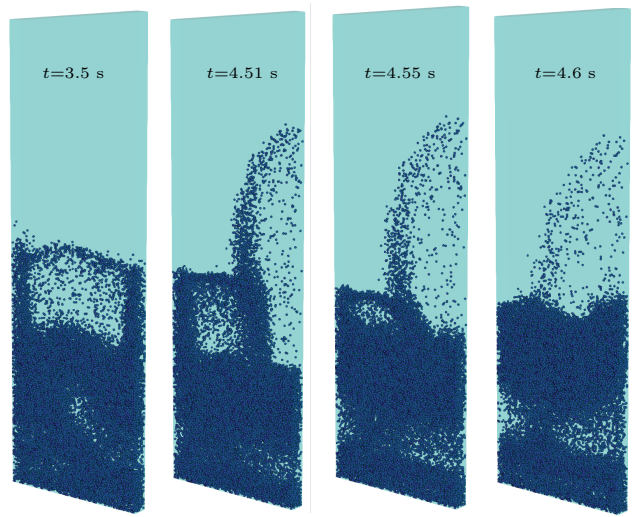
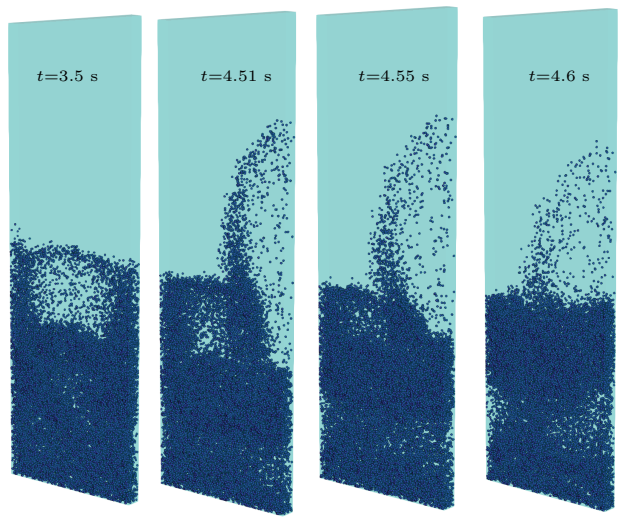
\begin{figure}[ht]
\vspace{1cm}
\begin{subfigure}{.55\textwidth}
  \begin{overpic}[width=\linewidth, grid=false]{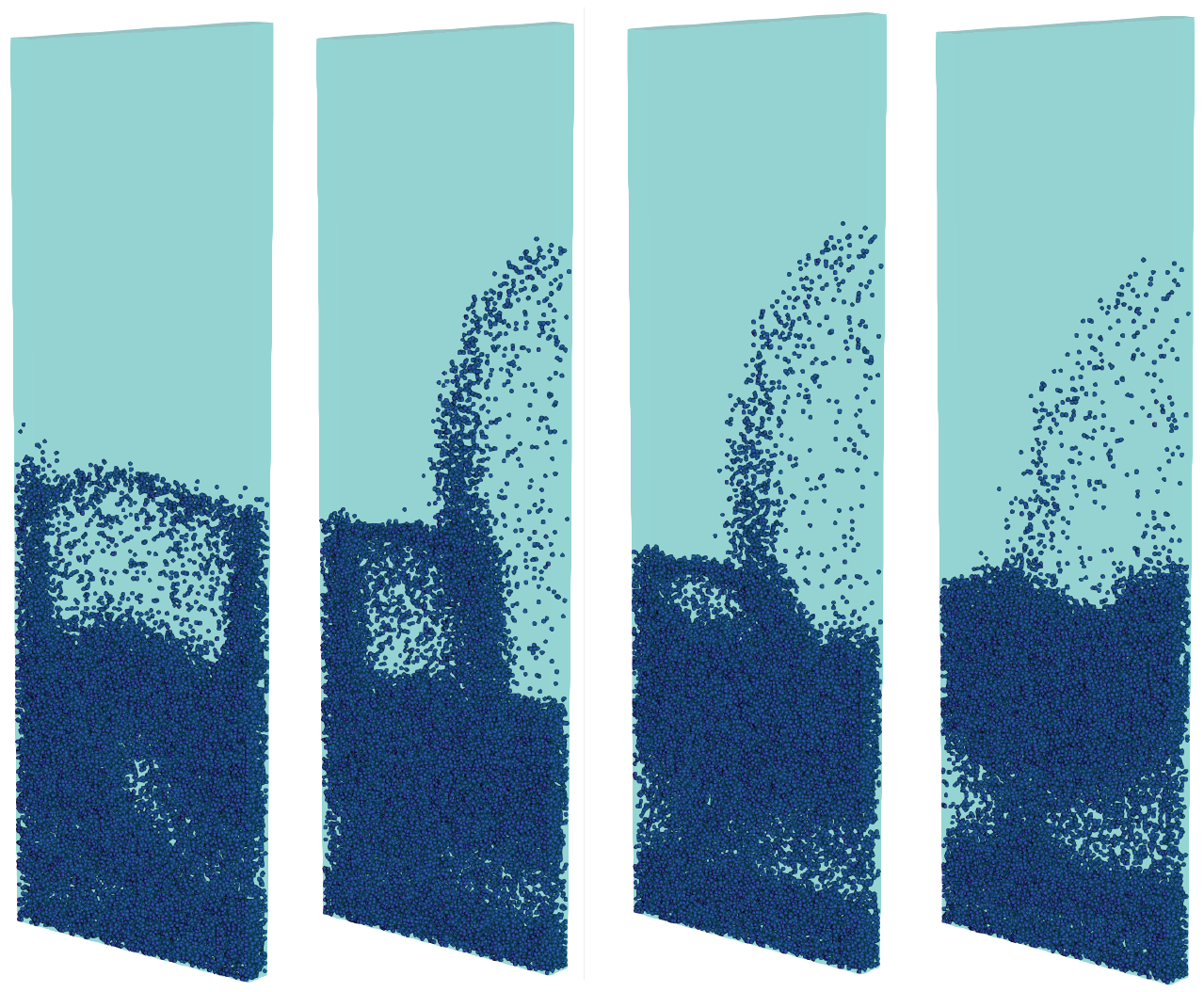}  
    \put(6,70){\tiny{\textcolor{black}{$t$=3.5 s}}}
    \put(32,70){\tiny{\textcolor{black}{$t$=4.51 s}}}
    \put(57,70){\tiny{\textcolor{black}{$t$=4.55 s}}}
    \put(84,70){\tiny{\textcolor{black}{$t$=4.6 s}}}
  \end{overpic}
  \caption{Unfiltered FOM solution at different time instances: $t=3.5$ s, $t=4.51$ s, $t=4.55$ s and $t=4.6$ s from left to right, respectively.}
\end{subfigure}
\par
\vspace{0.5cm}
\begin{subfigure}{.55\textwidth}
  \begin{overpic}[width=\linewidth, grid=false]{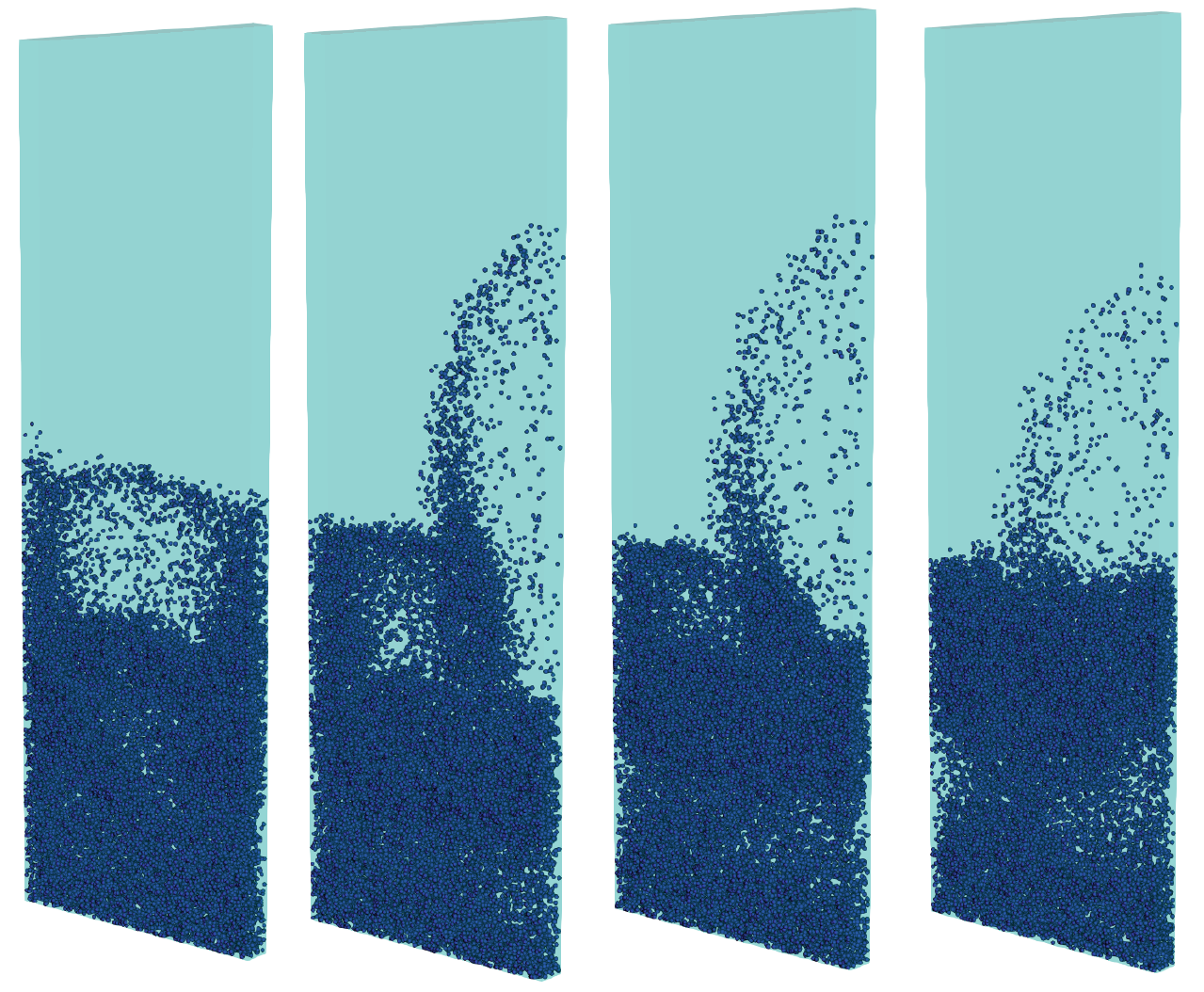}  
    \put(6,70){\tiny{\textcolor{black}{$t$=3.5 s}}}
    \put(32,70){\tiny{\textcolor{black}{$t$=4.51 s}}}
    \put(57,70){\tiny{\textcolor{black}{$t$=4.55 s}}}
    \put(84,70){\tiny{\textcolor{black}{$t$=4.6 s}}}
  \end{overpic}
  \caption{Filtered FOM solution at different time instances: $t=3.5$ s, $t=4.51$ s, $t=4.55$ s and $t=4.6$ s from left to right, respectively.}
\end{subfigure}

\caption{ROM validation - Lagrangian field: comparison of the time evolution of particle position computed by unfiltered (first row) and filtered (second row) FOM.}
\label{fig:Vis_FOM_denoised_Lagrangian}
\end{figure}

Fig. \ref{fig:PSD_lagrangian} shows the PSD as a function of frequency for two different particles, $l = 1000$ and $l = 24000$. The main figure is in logarithmic scale along $y$-axis. On the other hand, the inset is in logarithmic
scale for both $x$ and $y$ axes. The PSD threshold value, represented by the red horizontal line, is
set to $0.001$ for $\widetilde{y}$ and $\widetilde{z}$ components. For $\widetilde{x}$, we have considered a threshold of $5e-5$.  
Fig. \ref{fig:Vis_FOM_denoised_Lagrangian} shows a qualitative comparison between the original FOM snapshots and the filtered ones. We observe that
the agreement is good, so the selected PSD
threshold is well-suited for the system at hand.

\begin{figure}[ht]
\centering
\begin{subfigure}[b]{0.90\textwidth}
\centering
\subfloat[]{\label{fig:pred_z}\includegraphics[width=.4\linewidth]{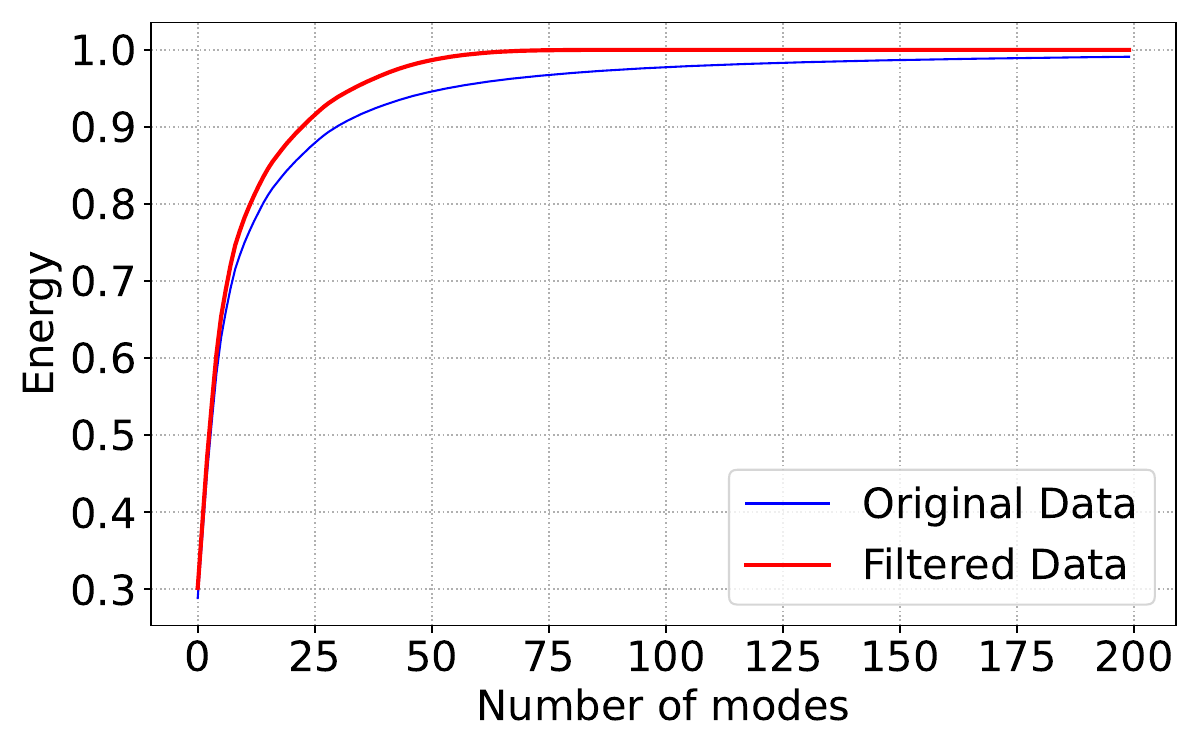}}
\subfloat[]{\label{fig:pred_y}\includegraphics[width=.4\linewidth]{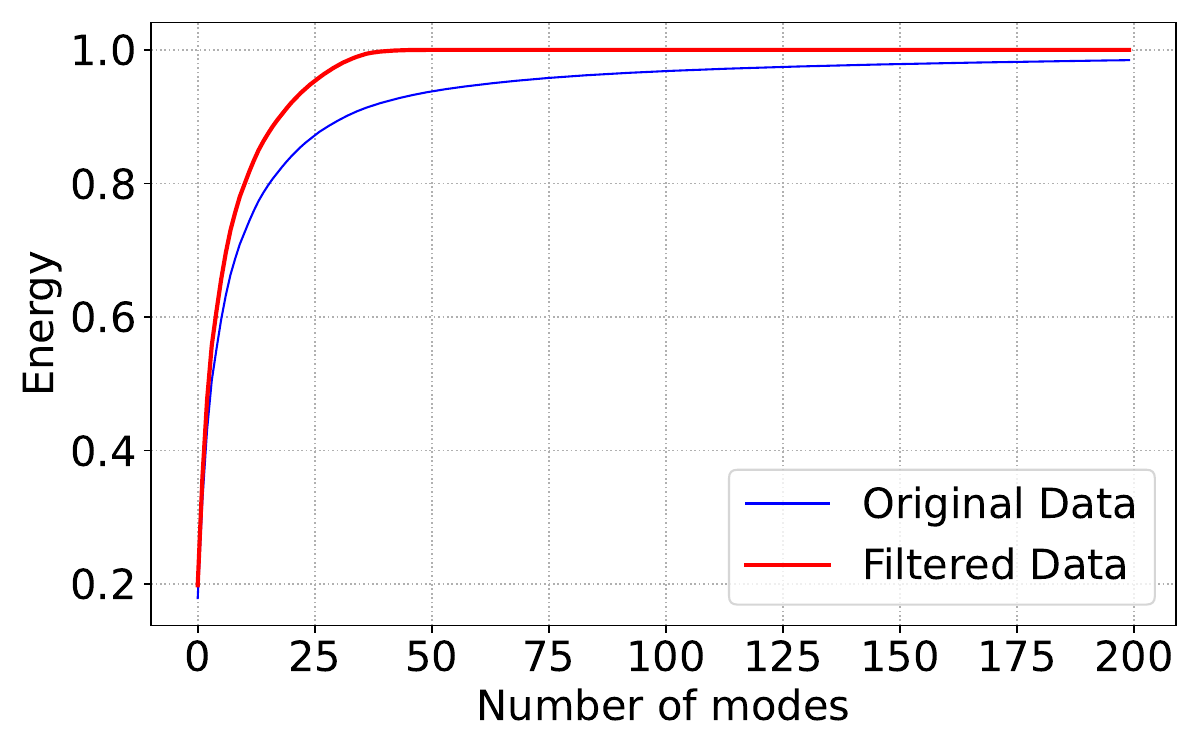}} \\ \vspace{0.5cm}
\subfloat[]{\label{fig:pred_x}\includegraphics[width=.4\linewidth]{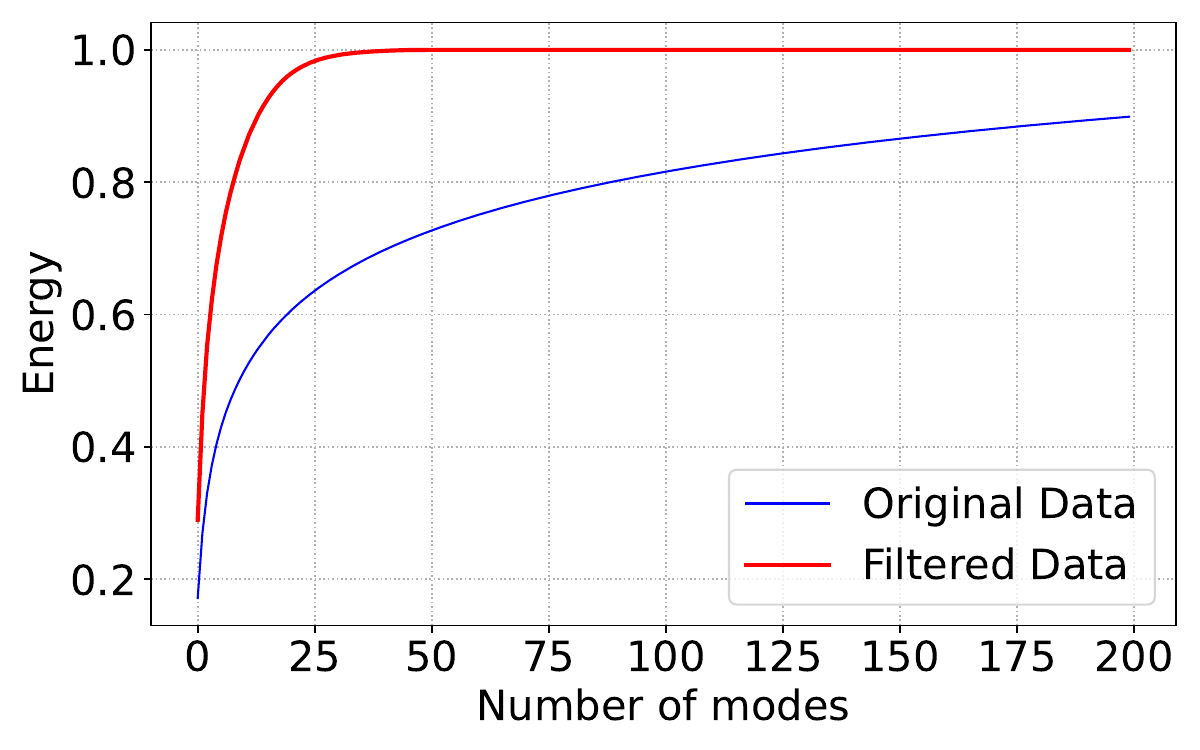}}
\end{subfigure}
\medskip
\caption{ROM validation - Lagrangian field: comparison of the cumulative energy $E$ defined in eq. \eqref{k_equation} of (a) $\widetilde{z}$, (b) $\widetilde{y}$ and (c) $\widetilde{x}$ for the original and filtered FOM snapshots.}
\label{fig:sigmaxyz}
\end{figure}

Fig. \ref{fig:sigmaxyz} shows the cumulative energy $E$ (eq. \eqref{k_equation}) for filtered and original FOM data. Like the Eulerian field, we note that the cumulative energy associated to the unfiltered snapshots shows a significantly slower convergence with respect to the filtered one.  

Table \ref{tab:Modes_num_L} reports
the number of modes associated to $\delta$ = 0.7, 0.9, 0.99: we can see that by performing the POD on
the original dataset, the cardinality of the latent space readuces drastically. This effect is stronger than the Eulerian field.  Infact, for $\delta = 0.99$, thanks to the filtering operation, we are able to reduce the number of modes of about 14-15 times for $\widetilde{x}$ and $\widetilde{y}$. For $\widetilde{x}$ the gain is evident also for the other two values of $\delta$, 0.9 and 0.7, where the number of modes is reduced of 10 and 5 times respectively. 
In order to maintain the trade-off between accuracy and efficiency, we intend to keep $25$ POD modes for $\widetilde{x}$ and $\widetilde{y}$ as well as $50$ POD modes corresponding to $\delta \approx 0.98$ for all the components. 

\begin{table}[ht!]
\centering
\begin{tabular}{|c|c|c|c|c|c|c|}
\hline  
   & Unfiltered $\widetilde{z}$  & Filtered $\widetilde{z}$ & Unfiltered $\widetilde{y}$  & Filtered $\widetilde{y}$ & Unfiltered $\widetilde{x}$  & Filtered $\widetilde{x}$ \\ \hline
$\delta = 70 \%$ & 9  &  8 & 10 & 8 & 42 & 6\\ \hline
$\delta = 90 \%$ & 31  & 24 &  33 & 19 & 202 & 14\\ \hline
$\delta = 99 \%$ & 195 &  54 & 202 & 14 & 452 & 30\\ \hline 
\end{tabular}
  \caption{ROM validation - Lagrangian field: number of required modes to retain different energy thresholds, $\delta = 0.7, 0.9, 0.99$, for original and filtered FOM data.
  }
 \label{tab:Modes_num_L}
\end{table}

Fig. \ref{fig:Vis_FOM_denoised_Lagrangian}  shows a visual FOM-ROM comparison for $t = 3.5$, 4.51, 4.55 and 4.6 s. For such time instances, the analysis is corroborated by the computation of the relative error 
\begin{equation}
E_{\widetilde{\bm{x}}}(t) = 100 \cdot \dfrac{||\widetilde{\bm{x}}_h(t) - \widetilde{\bm{x}}_r(t)||_{L^2(\Omega)}}{||{\widetilde{\bm{x}}_h}(t)||_{L^2(\Omega)}},
\label{eq:l2Error_lag}
\end{equation}
where $\widetilde{\bm{x}}_h$ is the particle position computed with the FOM and $\widetilde{\bm{x}}_r$ is the corresponding field computed with the ROM.
Note that $t = 3.5$ s belongs to the training set, so it is depicted to evaluate the ROM capability
to identify the system dynamics. The remaining three times, i.e., $t$ = 4.51, 4.55 and 4.6 s, are not associated with the training set and thus are used to check the accuracy of the ROM in predicting
the system dynamics. 
For $t = 3.5$ s we observe that our ROM approach is able to provide a very accurate reconstruction of the solution. Infact, we the $L^2$-norm error is of the $5\%$. Concerning the prediction phase, the overall comparison is still satisfactory  for $t = 4.51$ s and $t = 4.55$ with an error approximately of $8\%$ and $15\%$.  On the other hand, for $t= 4.6$ s, the solution provided by the surrogate model start to mismatch significantly with respect to the FOM one with an error reaching approximately $25\%$.  

\begin{figure}
\begin{subfigure}{.55\textwidth}
  \begin{overpic}[width=\linewidth, grid=false]{denoised2.pdf}  
    \put(6,70){\tiny{\textcolor{black}{$t$=3.5 s}}}
    \put(32,70){\tiny{\textcolor{black}{$t$=4.51 s}}}
    \put(57,70){\tiny{\textcolor{black}{$t$=4.55 s}}}
    \put(84,70){\tiny{\textcolor{black}{$t$=4.6 s}}}
  \end{overpic}
  \caption{FOM solution at different time instances: $t=3.5$ s, $t=4.51$ s, $t=4.55$ s and $t=4.6$ s from left to right, respectively.}
\end{subfigure}
\par
\vspace{0.5cm}
\begin{subfigure}{.55\textwidth}
 \begin{overpic}[width=\linewidth, grid=false]{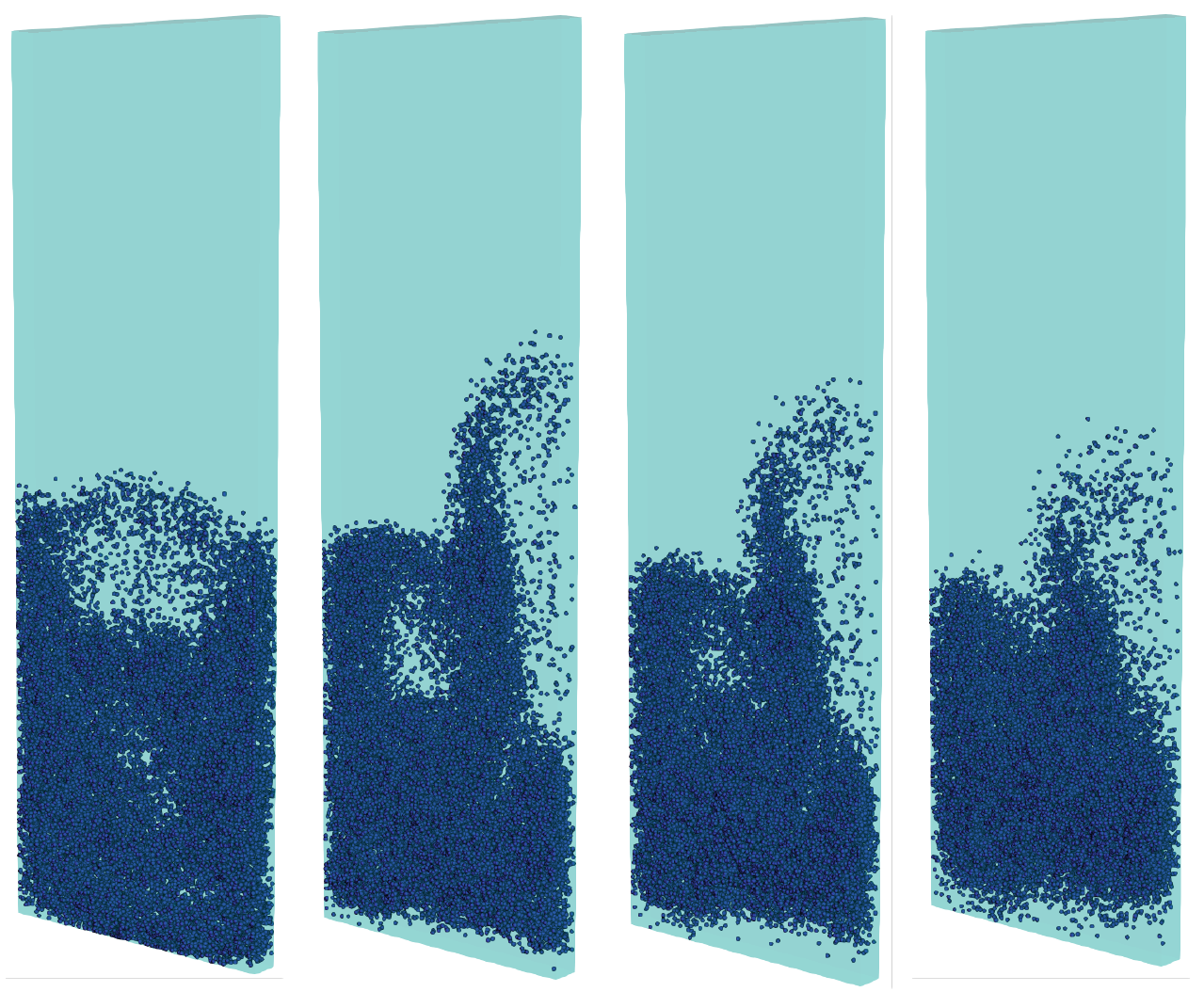}  
    \put(6,70){\tiny{\textcolor{black}{$t$=3.5 s}}}
    \put(32,70){\tiny{\textcolor{black}{$t$=4.51 s}}}
    \put(57,70){\tiny{\textcolor{black}{$t$=4.55 s}}}
    \put(84,70){\tiny{\textcolor{black}{$t$=4.6 s}}}
  \end{overpic}
  \caption{ROM solution at different time instances: $t=3.5$ s, $t=4.51$ s, $t=4.55$ s and $t=4.6$ s from left to right, respectively.}
\end{subfigure}

\caption{ROM validation - Lagrangian field: comparison of the time evolution of particle position computed by FOM (first row) and by ROM (second row) for $\delta \approx 98\%$. Notice that $t = 3.5$ s belongs to the training set while the other time instances are in the validation set.}
\label{fig:Vis_FOM_ROM_Lagrangian}
\end{figure}
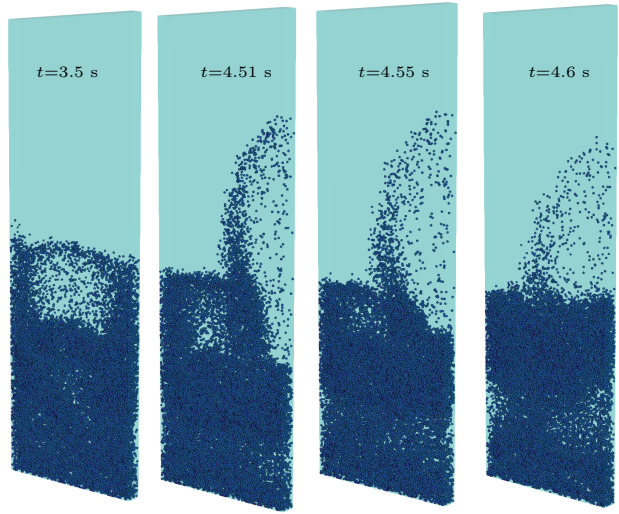
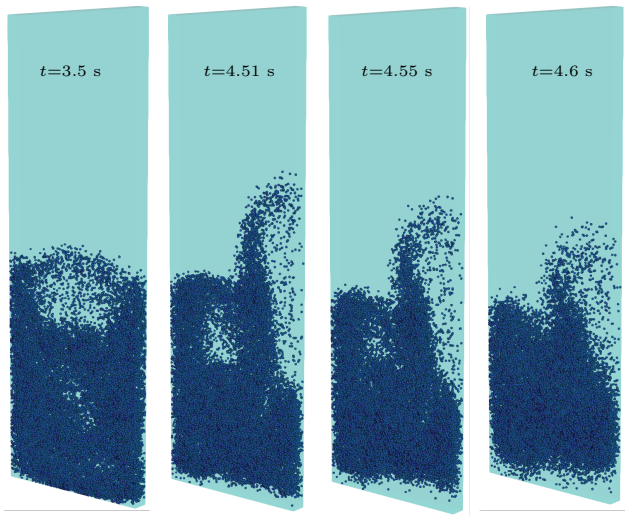

\subsection{Computational cost}\label{sec:cost}
We discuss briefly the efficiency of our ROM approach. We carried out FOM and ROM simulations on an 11th Gen Intel(R) Core(TM) i7-11700 @ 2.50GHz 32GB RAM by using one only processor. The CPU time associated with the FOM simulation is around 1.8e5 s, while the computation of reduced coefficients takes approximately 85 s for $\varepsilon$ and 113.58 s, 91.12 s and 101.01 s for $\widetilde{z}$,  $\widetilde{y}$ and $\widetilde{x}$, respectively. 
Therefore, the speed up, i.e. the ratio between the CPU time taken by the FOM simulation and the CPU time taken by the online phase, is of the order of 1e3 for all the variables.

\section{Concluding remarks}\label{sec:conc}

This work presents a non-intrusive data-driven ROM for fast and reliable CFD-DEM simulations.  
Since the development of ROMs in this framework is still in its infancy, our work in this paper has some novelty elements, among which the use of LSTM network for fluid-particle system prediction, the adoption of an FFT filtering approach to reduce the frequency content of the full order snapshots and a preliminary extension to the Lagrangian field. 

We assessed our ROM approach through a classical benchmark adopted for the validation of CFD-DEM solvers: a fluidized bed two-phase flow system. 
We found that for the Eulerian field our ROM is able to perform system prediction with a very good accuracy for short time windows. It can also capture the main unsteady flow features for longer time periods. On the other hand, for the Lagrangian field the results are not bad at all although they need to be improved. 
Moreover, our ROM results to be up to $1e-3$ times faster than the standard CFD-DEM approach. This makes the approach proposed very appealing both for and industrial applications where fluid-particle systems are involved.




As a follow-up of this work, we are going to improve the accuracy of our ROM approach for Lagrangian phase prediction. We think that some recent machine learning techniques, such as transformers \cite{jiang2023transcfd, wang2023swin}, could be adopted at this aim. 


\section{acknowledgments}\label{sec:akw}

We acknowledge the support provided by PON “Research and Innovation on Green related
issues” FSE REACT-EU 2021 project, INdAM-GNCS 2019-
2021 projects, SISSA-Dompè project “Development of a CFD model for INNOJET VENTILUS V1000 granulator”, PRIN “FaReX - Full and Reduced order modelling of coupled systems: focus on non-matching
methods and automatic learning” project, PNRR NGE iNEST “Interconnected Nord-Est Innovation
Ecosystem”.

\clearpage
\bibliographystyle{unsrt}
\bibliography{mybib}

\end{document}